\documentclass[10pt,a4paper]{article}
\makeatletter
\let\RerunFileCheck\@undefined  % 让它成为未定义状态
\let\RerunFileCheckOpen\@undefined
\makeatother
\usepackage[margin=3.24cm]{geometry}
\usepackage{graphicx} 
\usepackage{amsfonts,amsmath,amssymb,bm,dsfont,mathrsfs,physics}
\usepackage{amssymb,amsthm}
\usepackage{multirow}
\usepackage{makecell}
\usepackage{enumitem,cases}
\setlist[description]{style=unboxed,leftmargin=.5em}
\setlist[itemize]{style=sameline,leftmargin=2em}
\usepackage{float}
\usepackage[dvipsnames]{xcolor}
\colorlet{xlinkcolor}{red!50!black}

\definecolor{refblue}{RGB}{26,13,171}


\usepackage{hyperref}[6.83]
\hypersetup{
  colorlinks = true,
  allcolors = refblue,
  urlcolor = BrickRed,
}

\usepackage{cleveref}

\newtheorem{remark}{Remark}
\numberwithin{equation}{section}
\numberwithin{remark}{section}

\def\[{\begin{equation}}
\def\]{\end{equation}}
 
\setlist[description]{style=multiline,leftmargin=2em}
\setlist[itemize]{style=standard,leftmargin=2em}
\setlist[enumerate]{style=standard,leftmargin=2.2em,itemsep=3pt}

\usepackage{color}
\newcommand{\nn}{\nonumber}
\newlist{Steps}{enumerate}{2}
          \setlist[Steps]{label=\textbf{Step \arabic*.},itemindent=*,leftmargin=0pt}
\newcommand{\fall}{\quad\text{for all}\quad}
\newcommand{\re}{\mathrm{Re}}
\newcommand{\im}{\mathrm{Im}}

\usepackage{charter}
\usepackage[charter]{mathdesign}

\usepackage{algorithm}
\usepackage{algorithmic}

\title{\Large\bf
A Reduced Order Modeling Method with Variable-Separation-Based Domain Decomposition for Parametric Dynamical Systems \thanks{The research of this work was supported by 
the National Key R\&D Program of China (No. 2021YFA001300),
the National Natural Science Foundation of China (Nos. 12401567, 12271150, 12471405),
the Hunan Provincial Natural Science Foundation of China (No. 2023JJ10001), and 
the Science and Technology Innovation Program of Hunan Province (No. 2022RC1190).}}

\author{
Yuming Ba \thanks{School of Mathematics and Systems Science, Guangdong Polytechnic Normal University, Guangzhou 510665, China (\url{ yumingba@gpnu.edu.cn}).}
\and Liang Chen\thanks{School of Mathematics, Hunan University, Changsha 410082, China
(\url{chl@hnu.edu.cn}, \url{yrchen@hnu.edu.cn}, \url{qli28@hnu.edu.cn}).}
\thanks{Hunan Provincial Key Laboratory of Intelligent Information Processing and Applied Mathematics, Changsha 410082, China.}
\and Yaru Chen\footnotemark[3] \footnotemark[4]
\and Qiuqi Li\footnotemark[3] \footnotemark[4]
}

\begin{document}
\maketitle
%
%%%%% Begin Abstract %%%%%%%%%%%
\begin{abstract}
This paper proposes a model order reduction method for a class of parametric dynamical systems.
Using a temporal Fourier transform, we reformulate these systems into complex-valued elliptic equations in the frequency domain, containing frequency variables and parameters inherited from the original model. To reduce the computational cost of the frequency-variable elliptic equations, we extend the variable-separation-based domain decomposition method to the complex-valued context, resulting in an offline-online procedure for solving the parametric dynamical systems.
At the offline stage, separate representations of the solutions for the interface problem and the subproblems are constructed.
At the online stage, the solutions of the parametric dynamical systems for new parameter values can be directly derived by utilizing the separate representations and implementing the inverse Fourier transform. The proposed approach is capable of being highly efficient because the online stage is independent of the spatial discretization.
Finally, we present three specific instances of parametric dynamical systems to demonstrate the effectiveness of the proposed method.
\end{abstract}
%%%%% end %%%%%%%%%%%

%%%%% Keywords %%%%%%%%%%%

%%%% AMS subject classifications %%%%

%%%% maketitle %%%%%
\maketitle

%%%% Start %%%%%%

\section{Introduction}
\setcounter{equation}{0}
Parametric dynamical systems have been widely adopted in science and engineering to model complex real-world problems characterized by various uncertainties.
These uncertainties may stem from multiple sources, including physical properties, geometric configurations, initial conditions, and boundary conditions.
The numerical simulation of such parametric systems poses substantial computational challenges, especially in the case of repeated parameter evaluation for design optimization, control analysis, and uncertainty quantification.
Although traditional high-fidelity numerical approaches such as finite element method (FEM), finite difference method, and finite volume method can provide accurate solutions, they often entail extreme-scale computations with prohibitive computational costs.

To overcome arising computational bottlenecks in the simulation of parametric dynamical systems, a variety of innovative numerical approaches have been developed over past few decades. In particular, the model order reduction (MOR) technique has emerged as a powerful tool --- cf. Refs. \cite{Sapsis2009dynamical,Benner2015a,Brunton2019data,li2021non,R2019reduced,M2019data,Sleeman2022goal}.
It is used to construct reduced-order models by approximating solutions in a low-dimensional subspace, which significantly reduces the computational cost. A prominent MOR approach is the reduced basis (RB) method --- c.f. \cite{Haasdonk2008reduced,M2009A,Haasdonk2011efficient,Drohmann2012Reduced,hesthaven2016certified,hesthaven2022reduced}, which follows an offline-online computational strategy. In the offline phase, high-fidelity solutions for a carefully selected parameter set are precomputed and stored. The online phase then efficiently approximates solutions for new parameters via a linear combination of these precomputed solutions.
In this method, the key challenge lies in constructing an optimal set of basis functions that accurately capture the essential features of the full order model.
Having the reduced basis functions obtained, one can use intrusive or non-intrusive approach to construct a reduced order model. Non-intrusive MOR shows the flexibility and efficiency because the full order numerical systems are not required. To enhance the applicability of MOR to nonlinear systems and describe relationship between the inputs and outputs, MOR methods based on various machine or deep learning techniques have been proposed  \cite{Chen2021physical,Swischuk2019projection,Xiao2019error,kim2022a}. These methods roughly contain data-driven machine learning  \cite{Lee2019data,Mohebujjaman2019physically,Xu2020multi}, physical-informed machine learning  \cite{Chen2021physical,Swischuk2019projection,kim2022a}, and physics-data combined machine learning  \cite{Fu2023physics,Pan2024domain}. Among these methods, the most effective and easiest approach is the data-driven machine learning, which usually converges rapidly because the loss functions are often simple. The physically-informed machine learning can demonstrate strong theoretical constraints and inductive biases by integrating governing physical rules and domain knowledge into the learning process. The physics-data combined machine learning has been used to study nonlinear dynamical systems in small-data regimes \cite{Fu2023physics}. It has adopted a step-by-step training scheme, which seamlessly integrates the governing physical laws and the limited labeled data into feedforward neural networks.

An alternative approach for effective reducing the computational cost of dynamical systems is the frequency-domain method, which relies on two key transformations --- viz. the Fourier transform  \cite{Douglas1994approximation,Douglas1993frequency,Feng1994an,lee1999a,lee2000frequency} and the Laplace transform  \cite{lee2004an,lee2006a} in the time dimension.
The Fourier transform stands out as a significant discovery in mathematical sciences, playing a crucial role in modern scientific and technological advancements, which has been widely applied in the analysis of continuous-time systems, including signal and image processing, analog circuit analysis, and communication systems.
The fundamental principle behind this approach is to convert time-dependent PDEs into frequency-dependent elliptic problems, which are stationary in time but parameterized by frequency.
The time-domain solution is then reconstructed by applying an inverse Fourier transform to the frequency-domain solutions.
A major advantage of this framework is its inherent suitability for parallel computation --- i.e. each frequency can be solved independently as the resulting frequency-domain equations are decoupled.
For the frequency-domain equations that require numerous solutions, the MOR technique has demonstrated its effectiveness \cite{chen2010certified, Hess2013fast, Benner2015a}.
A low-order Galerkin proper orthogonal decomposition basis is employed to approximate solutions for frequency-domain equations in \cite{seo2023reduced}.
It adopts the Gaussian quadrature rule based on Legendre-Gauss-Lobatto (LGL) points to enable a precise numerical integration of the inverse Fourier transformation.

Motivated by these developments, we focus on the frequency-domain method based on the Fourier transform.
Specifically, we reformulate time-dependent parametric PDEs as a family of frequency-dependent parametric elliptic equations.
Although this transformation eliminates the time variable, the resulting equations remain dependent on spatial discretization, and must be solved repeatedly for multiple frequencies to recover the full solution via inverse transformation.
However, for large number of frequencies the overall computational cost is prohibitively high. Hence, the reduction of the model order is essential for improving efficiency.
As one of the MOR methods, the variable-separation (VS) method \cite{Li2017a} has been verified to be able to provide accurate and efficient approximate solutions to stochastic PDEs. Then, \cite{Jiang2018model} has extended the VS to stochastic saddle point problems. For nonlinear PDEs with random input, the VS method also very efficient \cite{Li2020a}. Moreover, it has been used in non-overlapping domain decomposition methods (DDM) for stochastic PDEs \cite{chen2024stochastic} and in parametric dynamical systems \cite{Chen2025a}.
To improve the computational efficiency of the frequency-variable parametric elliptic equations, we consider a variable-separation-based domain decomposition (DD-VS) method \cite{chen2024stochastic} instead of VS method. It is well-suited to heterogeneous problems and enables efficient online computations.

The DD-VS method is developed based on the framework of non-overlapping DDM \cite{Chan1987analysis, Sarker2009dd, waad2014pre, Mo2014inter, Liao2015ddua, Chen2015local, Xiao2019a, Xiao2019a2, Pan2024domain}.
In this approach, the computational domain is split into non-overlapping subdomains, where the degree of freedom in each subdomain are classified into interior and interface parts. Applying block Gaussian elimination, one reduces global system to a Schur complement system with interface unknowns only \cite{Mandel1993BalancingDD, Heinkenschloss2006neumann, J2001efficient}. This step constitutes the core and the most technically challenging part of the entire method.
The DD-VS method constructs a mapping between the input parameters and the interface system of the parametric PDEs. It adopts an offline-online decomposition strategy to improve computational efficiency.
In the offline phase, the extended VS method is used to transform the interface problem into a reduced-order algebraic system. A surrogate model is then built for this reduced interface system.
In the online phase, the surrogate model is utilized to rapidly evaluate new parameter samples. The online computation is highly efficient, as its cost is entirely independent of the spatial discretization.

To apply the DD-VS method to frequency-dependent elliptic equations, it is first necessary to extend the VS method to handle complex-valued problems, since the solutions of such equations, as well as their associated interface and subdomain problems, are inherently complex-valued.
To address this, we decompose the complex-valued problem into its real and imaginary parts, and formulate a coupled real-valued system accordingly. By following the standard procedures of the VS method developed for real-valued elliptic equations, we construct separate affine surrogate models for the real and imaginary components.
This extension allows the DD-VS method to be effectively applied to complex-valued frequency-dependent elliptic problems, while preserving computational efficiency in both the interface and subdomain problems during the online phase.

By combining the complex-valued DD-VS method with Fourier transformation, we develop a reduced order method for time-dependent problems, which can be regarded as an extension of the DD-VS method.
The process begins with applying the Fourier transform to reformulate the time-dependent problems into a set of frequency-variable elliptic equations.
Subsequently, utilizing the complex-valued DD-VS method, we construct efficient surrogate models for both the interface and subdomain problems during the offline phase.
The resulting online stage is highly efficient --- viz. for any new input parameter, the computation involves only evaluating the parametric coefficients of the surrogate model and performing the inverse Fourier transform. Importantly, the online phase is completely independent of the spatial discretization and requires no additional time-domain simulation, thereby enabling efficient solutions of time-dependent problems.

This paper is organized as follows. Section \ref{ssec:prelim} contains necessary notation and preliminaries and introduces the Fourier transform for the time variable.
Section \ref{sec-vs-ddm} introduces a VS method for parametric complex-valued systems and the domain decomposition method for deterministic complex-valued PDEs.
In Section \ref{sec-sddm}, we present a new DD-VS method with key steps and procedures.
Three numerical examples in Section \ref{sec-numerical examples} demonstrate the performance and computational advantages of the method proposed.
Finally, we draw conclusions and offer some remarks on the method and its potential applications.

\section{Preliminaries}
\label{ssec:prelim}
\setcounter{equation}{0}
Let $\Omega$ be the set representing the parameter domain. We use $\bm{\xi} \in \Omega$ to represent a vector of parametric inputs.
Given an open and bounded spatial domain $D\subset\mathbb{R}^{d}(d=1,2,3)$ with Lipschitz continuous boundary $\partial D$, we consider the parametric dynamical system
\begin{align*}
&\frac{\partial u}{\partial t}(\bm{x},t;\bm{\xi})=\mathcal{F}(u(\bm{x},t;\bm{\xi});\bm{\xi}), && \bm{x} \in D,~ t\in [0,T], ~ \bm{\xi} \in \Omega,\\
&u(\bm{x},0;\bm{\xi})= u_0(\bm{x};\bm{\xi}), &&  \bm{x} \in D, ~ \bm{\xi} \in \Omega,\\
&\mathcal{B}(u(\bm{x},t;\bm{\xi}))=g(\bm{x},t;\bm{\xi}), && \bm{x} \in \partial D, ~ t\in [0,T], ~ \bm{\xi} \in \Omega,
\end{align*}
where $T$ is the end time of the computation, $\mathcal{F}$ the spatial differential operator, $\mathcal{B}$ the boundary condition operator,
$u(\cdot,t;\bm{\xi})$ the system solution, $g(\cdot,t;\bm{\xi})$ the boundary term, and $u_0(\cdot;\bm{\xi})$ the initial function.
To simplify the presentation, we illustrate our methodology using the following parametric parabolic partial differential equation as a representative example
\begin{equation}\label{eq-dynamical-system}
\begin{aligned}
   & \frac{\partial u}{\partial t}(\bm{x},t;\bm{\xi})-\nabla\cdot (c(\bm{x};\bm{\xi}) \nabla u(\bm{x},t;\bm{\xi}))=f(\bm{x},t;\bm{\xi}), &&\bm{x} \in D,~ t\in [0,T], ~ \bm{\xi} \in \Omega,\\
   & u(\bm{x},0;\bm{\xi})=u_0(\bm{x};\bm{\xi}), &&  \bm{x} \in D, ~ \bm{\xi} \in \Omega,\\
    &u(\bm{x},t;\bm{\xi})=g(\bm{x},t;\bm{\xi}), && \bm{x} \in \partial D,~ t\in [0,T], ~ \bm{\xi} \in \Omega,
\end{aligned}
\end{equation}
where $c(\bm{x};\bm{\xi})$ is the diffusion coefficient and $f(\bm{x},t;\bm{\xi})$ the source term.

Let
\begin{equation*}
\mathcal{V}:= L^2(D;\mathbb{C}) = \{\hat{u}(\bm{x})|\int_D |\hat{u}(\bm{x})|^2 dx < + \infty \}.
\end{equation*}
Besides, the inner product on $\mathcal{V}$ and the associated norm are defined by
\begin{equation*}
\langle \hat{u}(\bm{x}),
\hat{v}(\bm{x}) \rangle := \int_D \hat{u}(\bm{x}) \overline{\hat{v}(\bm{x})} d\bm{x},
\quad
\|\hat{u}(\bm{x})\|:=\sqrt{\langle \hat{u}(\bm{x}),\hat{u}(\bm{x})\rangle}.
\end{equation*}

\subsection{The frequency-domain method}
To tackle the dynamical system \eqref{eq-dynamical-system}, we apply the Fourier transformation to get the following set of complex-valued elliptic equations depending on the frequency $\omega$
\begin{equation}\label{eq-fourier}
\begin{aligned}
&i\omega \hat{u}(\bm{x},\omega;\bm{\xi}) - \nabla \cdot (c(\bm{x};\bm{\xi})\nabla \hat{u}(\bm{x},\omega;\bm{\xi})) = \hat{f}(\bm{x},\omega;\bm{\xi}),  &&  \bm{x} \in D, ~ \bm{\xi} \in \Omega,\\
&\hat{u}(\bm{x},\omega;\bm{\xi}) = \hat{g}(\bm{x},\omega;\bm{\xi}), &&  \bm{x} \in \partial D,~ \bm{\xi} \in \Omega,  
\end{aligned}
\end{equation}
where $i=\sqrt{-1}$ is the imaginary unit. The Fourier transform $\hat{u}(\bm{x},\omega;\bm{\xi})$ of a function $u(\bm{x},t;\bm{\xi})$ in time and its inverse are
\begin{align*}
\hat{u}(\bm{x},\omega;\bm{\xi}) &:= \int_{-\infty}^{\infty} u(\bm{x},t;\bm{\xi}) \exp(-i\omega t) \mathrm{d} t,
\\[1mm]  u(\bm{x},t;\bm{\xi}) &:= \frac{1}{2\pi} \int_{-\infty}^{\infty} \hat{u}(\bm{x},\omega;\bm{\xi}) \exp(i\omega t) \mathrm{d} \omega.
\end{align*}
For the Fourier transformation, the function $u$ is extended by zero outside of the interval $[0, T]$ with the functions $f$ and $g$ being treated identically to derive equation \eqref{eq-fourier}.

The weak formulation of \eqref{eq-fourier} is: For all $\omega\in\mathbb{R}$ and $\bm{\xi}\in\Omega $, find $\hat{u}(\cdot,\omega;\bm{\xi})$ such that  $\hat{u} = \hat{g}$ on $\partial D$ and
\begin{equation}\label{eq-fourier-weak}
 i\omega \langle \hat{u}(\cdot,\omega;\bm{\xi}), \hat{v}\rangle
 + \langle c(\cdot;\bm{\xi})\nabla \hat{u}(\cdot,\omega;\bm{\xi}), \nabla \hat{v}\rangle =  \langle  \hat{f}(\cdot,\omega;\bm{\xi}), \hat{v} \rangle
 \fall \hat{v}\in\mathcal{V}.
\end{equation}

Consider the finite element approximation of problem \eqref{eq-fourier-weak} in an $n$-dimensional subspace $\mathcal{V}_h\subset \mathcal{V}$. If $\{\psi_{j}\}_{j=1}^{n}$ is a set of basis functions of the space $\mathcal{V}_h$, then the solution $\hat{u}(\bm{x},\omega;\bm{\xi})$ can be approximated as
\begin{equation*}
    \hat{u}(\bm{x},\omega;\bm{\xi})\approx \hat{u}_h(\bm{x},\omega;\bm{\xi}):=\sum_{j=1}^{n}\hat{u}_{j}(\omega;\bm{\xi}) \psi_{j}(\bm{x}),
\end{equation*}
where 
\begin{equation*}
\hat{u}_{j}(\omega;\bm{\xi}):=\hat{u}_{j}^{\re}(\omega;\bm{\xi}) +i \hat{u}_{j}^{\im}(\omega;\bm{\xi}),\quad j=1, \ldots, n,
\end{equation*}
is a complex-valued function.
For the approximation of the  Fourier inverse,
we employ the Gaussian quadrature  based on the LGL-points with an appropriate interval $[0, \omega^*]$ for a sufficiently large $\omega^*$, where $\hat{u}(\cdot,\omega;\cdot)$ is negligible for $|\omega|>\omega^*$. Let $\{\omega_j\}_{j=1}^{N_\omega}$ be the set of LGL-points on the interval $[0,\omega^*]$. Then, the time-variable approximate $u(\bm{x},t;\bm{\xi})$ to problem \eqref{eq-dynamical-system} is
\begin{equation}\label{eq:ift}
    \quad u(\bm{x},t;\bm{\xi}) = \frac{1}{\pi} \re\bigg(\sum_{j=1}^{N_\omega} \hat{u}_h(\bm{x},\omega_j;\bm{\xi}) \exp(i\omega_j t) w_j\bigg),
\end{equation}
where $w_j$ is the Gaussian quadrature weight of the LGL-points $\omega_j$.

Note that, for different values of $\omega$ and $\bm{\xi}$, one has to solve \eqref{eq-fourier-weak} repetitively to obtain the solution $\hat{u}_h(\bm{x},\omega;\bm{\xi})$. This process requires a large amount of computation when dealing with a large number of parameters $\omega$ and $\bm{\xi}$, especially for more complex problems.
To reduce the computational cost, the frequency variable $\omega$ will be treated as a one-dimensional parameter in this paper.
We aim to construct the reduced order model of $\hat{u}_h(\bm{x},\omega;\bm{\xi})$.
With this reduced order model, \eqref{eq:ift} can be computed in a computationally economical fashion, and the implementation details will be elaborated in the following sections.

\section{VS and Domain Decomposition Methods}
\label{sec-vs-ddm}
In this section, we will present in detail the VS method and the domain decomposition method for the complex-valued problem \eqref{eq-fourier-weak}, where the VS method will be utilized in the next section to construct surrogate models for domain decomposed problems.

\subsection{VS method for complex-valued problems}
\label{sec-vs}
As already mentioned, our goal is the construction of a reduced order model when $\omega$ is considered as one of the parameters.
In this section, the VS method will be extended to derive a separate representation of the solutions $\hat{u}(\bm{x},\omega;\bm{\xi})$ for the complex-valued elliptic problem \eqref{eq-fourier-weak}. For convenience, define $\bm{\mu}:=(\omega,\bm{\xi})$, $\tilde\Omega:=\mathbb{R}\times\Omega$. Then the equality \eqref{eq-fourier-weak} can be rewritten as
\begin{equation}\label{eq-fourier-weak-2}
a(\hat{u},\hat{v};\bm{\mu}) = b(\hat{v};\bm{\mu}),
\end{equation}
where $a(\cdot,\cdot;\bm{\mu})$ and $b(\cdot;\bm{\mu})$ are complex-valued bilinear and linear forms on $\mathcal{V}$, respectively.
Assume that $a(\cdot,\cdot;\bm{\mu})$ and $b(\cdot;\bm{\mu})$ are affine with respect to $\bm{\mu}$, i.e.
\begin{equation}
\begin{aligned}\label{eq-fourier-affine}
a(\hat{u},\hat{v};\bm{\mu})&=\sum_{j=1}^{m_a}\alpha_{j}(\bm{\mu})a_{j}(\hat{u},\hat{v})+i\gamma(\bm{\mu})m(\hat{u},\hat{v}) ~\ \fall \hat{u},\hat{v}\in\mathcal{V},  \ \bm{\mu}\in \tilde{\Omega},\\
b(\hat{v};\bm{\mu})&=\sum_{j=1}^{m_b}\bigg(\beta_{j}^{\re}(\bm{\mu})b_{j}^{\re}(\hat{v})+i\beta_{j}^{\im}(\bm{\mu})b_{j}^{\im}(\hat{v})\bigg)~\ \fall \hat{v}\in\mathcal{V}, \ \bm{\mu}\in \tilde{\Omega},
\end{aligned}
\end{equation}
where $\gamma(\bm{\mu}) = \omega$, $m(\hat{u},\hat{v})=\langle \hat{u},\hat{v} \rangle$, $\alpha_{j}(\bm{\mu}):\tilde\Omega\to R$ is $\bm{\mu}$-dependent function and  $a_{j}:\mathcal{V}\times\mathcal{V}\to R$ is a bilinear form independent of $\bm{\mu}$, for each $j=1,\ldots,m_a$, $\beta_{j}^{\re}(\bm{\mu}):\tilde\Omega\to R$ and $\beta_{j}^{\im}(\bm{\mu}):\tilde\Omega\to R$ are $\bm{\mu}$-dependent functions, $b_{j}^{\re}:\mathcal{V}\to R$ and $b_{j}^{\im}:\mathcal{V}\to R$, $j=1,\ldots,m_b$are linear forms independent of $\bm{\mu}$.
Our analysis shows that affine decompositions play a crucial role in enabling the offline-online decomposition.

When \eqref{eq-fourier-affine} holds, the matrix form of Eq.~\eqref{eq-fourier-weak} in a finite element space $\mathcal{V}_h$ is given by
\begin{equation}\label{eq-matrix-equations}
    \left(\sum_{j=1}^{m_{a}}\alpha_{j}(\bm{\mu})A_{j} + i \gamma(\bm{\mu}) M \right)\hat{\mathbf{u}}(\bm{\mu})=\sum_{j=1}^{m_b}\left(\beta_{j}^{\re}(\bm{\mu})F^{\re}_j+i\beta_{j}^{\im}(\bm{\mu})F^{\im}_j\right),
\end{equation}
where $\hat{\mathbf{u}}(\bm{\mu}) = \hat{\mathbf{u}}^{\re}(\bm{\mu})+i\hat{\mathbf{u}}^{\im}(\bm{\mu})$, and
\begin{align*}
 &(A_{j})_{kl}=a_{j}(\psi_{k},\psi_{l}), \quad M_{kl}=\langle\psi_{k},\psi_{l}\rangle, \\ & (\hat{\mathbf{u}}^{\re}(\bm{\mu}))_{k}=\hat{u}_{k}^{\re}(\bm{\mu}), \quad (\hat{\mathbf{u}}^{\im}(\bm{\mu}))_{k}=\hat{u}_{k}^{\im}(\bm{\mu}),\\
&(F_{j}^{\re})_{k}=b_{j}^{\re}(\psi_{k}),\quad (F_{j}^{\im})_{k}=b_{j}^{\im}(\psi_{k}),\quad  1\leq k,l\leq n.
\end{align*}

Now, we introduce the VS method for the problem \eqref{eq-matrix-equations} to obtain a separate approximation of the solution --- i.e.
\begin{equation}\label{eq-approx}
\hat{\mathbf{u}}(\bm{\mu})\approx \hat{\mathbf{u}}_{N}(\bm{\mu}):=\sum_{j=1}^{N}\left(\zeta_j^{\re}(\bm{\mu}) \mathbf{c}_j^{\re}+i\zeta_j^{\im}(\bm{\mu}) \mathbf{c}_j^{\im}\right),
\end{equation}
where $\zeta_j^{\re}(\bm{\mu})$ and $\zeta_j^{\im}(\bm{\mu})$ are real-valued parametric functions, and $\mathbf{c}_j^{\re}$ and $\mathbf{c}_j^{\im}$ are real-valued vectors independent of $\bm{\mu}$, with $N(\ll n)$ denoting the number of the separate terms.
We define the error of the VS method by
\begin{equation}\label{eq-e}
\mathbf{e}(\bm{\mu}):=\mathbf{e}^{\re}(\bm{\mu}) + i\mathbf{e}^{\im}(\bm{\mu})=\hat{\mathbf{u}}(\bm{\mu})-\hat{\mathbf{u}}_{k-1}(\bm{\mu}),
\end{equation}
where $\hat{\mathbf{u}}_{k-1}(\bm{\mu}) \equiv 0$ at $k=1$.
Substituting \eqref{eq-e} into the Eq.~\eqref{eq-matrix-equations} yields
\begin{align}\label{eq-g-re}
- \gamma(\bm{\mu}) M \mathbf{e}^{\im}(\bm{\mu}) + \sum_{j=1}^{m_{a}}\alpha_{j}(\bm{\mu})A_{j}\mathbf{e}^{\re}(\bm{\mu}) &= \mathbf{r}_k^{\re}(\bm{\mu}), \\
\label{eq-g-im}
\gamma(\bm{\mu}) M \mathbf{e}^{\re}(\bm{\mu}) + \sum_{j=1}^{m_{a}}\alpha_{j}(\bm{\mu})A_{j} \mathbf{e}^{\im}(\bm{\mu}) &= \mathbf{r}_k^{\im}(\bm{\mu}),
\end{align}
where
\begin{equation}\label{eq-r-re}
\mathbf{r}_k^{\re}(\bm{\mu}):=
\left\{
\begin{aligned}
    &\sum_{j=1}^{m_b}\beta_{j}^{\re}(\bm{\mu})F^{\re}_j,&k=1,\\
    &\sum_{j=1}^{m_b}\beta_{j}^{\re}(\bm{\mu})F^{\re}_j + \gamma(\bm{\mu}) M  \hat{\mathbf{u}}_{k-1}^{\im}(\bm{\mu}) -
    \sum_{j=1}^{m_{a}}\alpha_{j}(\bm{\mu})A_{j}\hat{\mathbf{u}}_{k-1}^{\re}(\bm{\mu}),&k\geq2,
\end{aligned}
\right.
\end{equation}
and
\begin{equation}\label{eq-r-im}
\mathbf{r}_k^{\im}(\bm{\mu}):=
\left\{
\begin{aligned}
    &\sum_{j=1}^{m_b}\beta_{j}^{\im}(\bm{\mu})F^{\im}_j,&k=1,\\
    &\sum_{j=1}^{m_b}\beta_{j}^{\im}(\bm{\mu})F^{\im}_j- \gamma(\bm{\mu}) M \hat{\mathbf{u}}_{k-1}^{\re}(\bm{\mu})
    - \sum_{j=1}^{m_{a}}\alpha_{j}(\bm{\mu})A_{j}\hat{\mathbf{u}}_{k-1}^{\im}(\bm{\mu}),&k\geq2.
\end{aligned}
\right.
\end{equation}
At the $k$-th step  we choose $\bm{\mu}_k$ as follows:
\begin{equation*}
\bm{\mu}_k=
\left\{
\begin{aligned}
&\text{chosen randomly in } \Xi,&k=1,\\
&\mathop{\text{argmax}}\limits_{\bm{\mu}\in\Xi}|\mathbf{r}_k(\bm{\mu})|,&k\geq2,
\end{aligned}
\right.
\end{equation*}
where $\Xi$ is a collection of a finite number of samples in $\tilde{\Omega}$ and $\mathbf{r}_k(\bm{\mu}) = \mathbf{r}_k^{\re}(\bm{\mu}) + i\mathbf{r}_k^{\im}(\bm{\mu})$.
Solving the Eqs. \eqref{eq-g-re}-\eqref{eq-g-im} with $\bm{\mu} = \bm{\mu}_k$, we derive $\mathbf{c}_k^{\re}= \mathbf{e}^{\re}(\bm{\mu}_k) $ and $\mathbf{c}_k^{\im}= \mathbf{e}^{\im}(\bm{\mu}_k)$ in \eqref{eq-approx}.

Let $\tilde{\mathbf{e}}^{\re}(\bm{\mu}) = \zeta_k^{\re}(\bm{\mu}) \mathbf{c}_k^{\re}$ and $\tilde{\mathbf{e}}^{\im}(\bm{\mu}) = \zeta_k^{\im}(\bm{\mu}) \mathbf{c}_k^{\im}$.
By taking $\mathbf{e}^{\re} = \tilde{\mathbf{e}}^{\re}$ and $\mathbf{e}^{\im} = \tilde{\mathbf{e}}^{\im}$ in \eqref{eq-g-re}-\eqref{eq-g-im}, the linear systems involving the unknowns $\zeta_k^{\re}(\bm{\mu})$ and $\zeta_k^{\im}(\bm{\mu})$ are as follows
\begin{align}
\label{eq-xi-unknown-1}
- \gamma(\bm{\mu}) M \mathbf{c}_k^{\im}  \zeta_k^{\im}(\bm{\mu}) + \sum_{j=1}^{m_{a}}\alpha_{j}(\bm{\mu})A_{j} \mathbf{c}_k^{\re} \zeta_k^{\re}(\bm{\mu})&= \mathbf{r}_k^{\re}(\bm{\mu}), \\
\label{eq-xi-unknown-2}
\gamma(\bm{\mu}) M  \mathbf{c}_k^{\re} \zeta_k^{\re}(\bm{\mu}) + \sum_{j=1}^{m_{a}}\alpha_{j}(\bm{\mu})A_{j} \mathbf{c}_k^{\im} \zeta_k^{\im}(\bm{\mu}) &= \mathbf{r}_k^{\im}(\bm{\mu}).
\end{align}
Taking the dot product of \eqref{eq-xi-unknown-1} with $\mathbf{c}^{\re}_k$ and \eqref{eq-xi-unknown-2} with $\mathbf{c}^{\im}_k$ yields the affine representations of $\zeta_k^{\re}(\bm{\mu})$ and $\zeta_k^{\im}(\bm{\mu})$ as
\begin{align}
\label{eq-xi-1}
 &\quad \zeta_k^{\re}(\bm{\mu})\nn\\
 &= \frac{\sum_{j=1}^{m_{a}}\alpha_{j}(\bm{\mu}) (\mathbf{c}_k^{\im})^TA_{j} \mathbf{c}_k^{\im} \cdot (\mathbf{c}_k^{\re})^T\mathbf{r}_k^{\re}(\bm{\mu})+\gamma(\bm{\mu})(\mathbf{c}_k^{\re})^T M \mathbf{c}_k^{\im}\cdot (\mathbf{c}_k^{\im})^T\mathbf{r}_k^{\im}(\bm{\mu})}{\sum_{j=1}^{m_{a}}\alpha_{j}(\bm{\mu})(\mathbf{c}_k^{\re})^TA_{j} \mathbf{c}_k^{\re} \cdot \sum_{j=1}^{m_{a}}\alpha_{j}(\bm{\mu})(\mathbf{c}_k^{\im})^TA_{j} \mathbf{c}_k^{\im}+ \gamma(\bm{\mu})(\mathbf{c}_k^{\re})^T M \mathbf{c}_k^{\im} \cdot \gamma(\bm{\mu}) (\mathbf{c}_k^{\im})^T  M  \mathbf{c}_k^{\re} },\\
\label{eq-xi-2}
 &\quad \zeta_k^{\im}(\bm{\mu})\nn\\
 &=  \frac{\sum_{j=1}^{m_{a}}\alpha_{j}(\bm{\mu})(\mathbf{c}_k^{\re})^TA_{j} \mathbf{c}_k^{\re}\cdot (\mathbf{c}_k^{\im})^T \mathbf{r}_k^{\im}(\bm{\mu})-\gamma(\bm{\mu})(\mathbf{c}_k^{\im})^T M  \mathbf{c}_k^{\re} \cdot (\mathbf{c}_k^{\re})^T\mathbf{r}_k^{\re}(\bm{\mu}) }{\gamma(\bm{\mu})(\mathbf{c}_k^{\re})^T M \mathbf{c}_k^{\im} \cdot \gamma(\bm{\mu})(\mathbf{c}_k^{\im})^T M  \mathbf{c}_k^{\re} +\sum_{j=1}^{m_{a}}\alpha_{j}(\bm{\mu}) (\mathbf{c}_k^{\im})^T A_{j} \mathbf{c}_k^{\im} \cdot \sum_{j=1}^{m_{a}}\alpha_{j}(\bm{\mu}) (\mathbf{c}_k^{\re})^TA_{j} \mathbf{c}_k^{\re}},
\end{align}
where both $\mathbf{r}_k^{\re}(\bm{\mu})$ and $\mathbf{r}_k^{\im}(\bm{\mu})$ are affine functions of $\bm{\mu}$, as defined by the Eqs.~\eqref{eq-r-re}-\eqref{eq-r-im}.

The iteration procedure ends when $|\mathbf{r}_k(\bm{\mu}_k)|$ is small enough.
The above procedure of the VS method for the complex-valued elliptic problem is summarized in Algorithm \ref{algorithm-vs-sSAS}.

\begin{algorithm}
\caption{VS method for the complex-valued elliptic problem}
\label{algorithm-vs-sSAS}
\begin{algorithmic}[1]
\REQUIRE{The complex-valued elliptic problem \eqref{eq-matrix-equations}, a set of samples $\Xi\in \tilde\Omega$, and the error tolerance $\varepsilon$.}
\ENSURE{The approximation $\hat{\mathbf{u}}_{N}(\bm{\mu}):=\sum_{j=1}^{N}\Big(\zeta_j^{\re}(\bm{\mu}) \mathbf{c}_j^{\re}+i\zeta_j^{\im}(\bm{\mu}) \mathbf{c}_j^{\im}\Big)$.}
\begin{enumerate}
\item Initialize the iteration counter $k=1$, a random ${\bm{\mu}_1}\in \Xi$;
\item Calculate $\mathbf{c}_k^{\re}$ and $\mathbf{c}_k^{\im}$ by solving \eqref{eq-g-re} and \eqref{eq-g-im} with $\bm{\mu}={\bm{\mu}_k}$; compute $\zeta_k^{\re}(\bm{\mu})$ and $\zeta_k^{\im}(\bm{\mu})$ by \eqref{eq-xi-1} and \eqref{eq-xi-2};
\item Update $\Xi=\Xi \textbf{\textbackslash} \bm{\mu}_k $, and take $\hat{\mathbf{u}}_{k}(\bm{\mu}):=\sum_{j=1}^{k}\Big(\zeta_j^{\re}(\bm{\mu}) \mathbf{c}_j^{\re}+i\zeta_j^{\im}(\bm{\mu}) \mathbf{c}_j^{\im}\Big)$;
\item Set $k\to k+1$;\\
\item If $|\mathbf{r}(\bm{\mu})| \geq \varepsilon$ and $\Xi \neq \emptyset$, choose $\bm{\mu}_k\in \mathop{\text{argmax}}\limits_{\bm{\mu}\in\Xi}|\mathbf{r}_k(\bm{\mu})|$, and \textbf{return} to Step 2; otherwise set $N=k$ and \textbf{terminate}.
\end{enumerate}
\end{algorithmic}
\end{algorithm}

We emphasize that the coefficients $\{\mathbf{c}_j^{\re}\}_{j=1}^{N}$ and $\{\mathbf{c}_j^{\im}\}_{j=1}^{N}$ (independent of $\bm{\mu}$) can be calculated during the offline stage of the VS method and used directly
in the online stage. Then we only need to calculate the parametric functions in Eq.~\eqref{eq-approx} for any $\bm{\mu} \in \tilde\Omega$ during the online stage. This procedure is highly efficient, since it only involves the separate representation  \eqref{eq-approx}.

Although the VS method allows constructing reduced-order models for frequency-domain problems, its application to time-domain dynamical problems in complex geometries remains challenging. In particular, employing the Fourier transform in such scenarios often necessitates a large $N$ to achieve sufficient approximation accuracy, leading to high computational costs. To address this limitation and enhance computational efficiency, we propose incorporating a domain decomposition strategy.

\subsection{Domain decomposition method for complex-valued problems}\label{sec-ddm}

In this section, we  introduce the non-overlapping domain decomposition method for complex-valued problems. For real-valued problems, the reader can consult Refs.~\cite{waad2014dd, Mu2018ddrd}.
Here the complex-valued problem \eqref{eq-fourier} is considered for a fixed $\bar{\bm{\mu}}\in\tilde\Omega$.

Suppose that the domain $D$ is split into $N_s$ non-overlapping subdomains $\{D_j\}_{j=1}^{N_s}$, i.e. $D=\mathop{\large\cup}\limits_{j=1}^{N_s}\overline{D_j}$ and $D_j\cap D_k=\emptyset$ if $j\neq k$. Restricting the weak formulation \eqref{eq-fourier-weak-2} to subdomain $D_j$, we write
\begin{equation*}
    a(\hat{u},\hat{v};\bar{\bm{\mu}})_j = b(\hat{v};\bar{\bm{\mu}})_j~\ \fall \hat{v}\in \mathcal{V}_j,
\end{equation*}
where $\mathcal{V}_{j}$ is the Hilbert space $\mathcal{V}$ restricted to the subdomain $D_j$.
The finite element discretization of the above equation leads the local linear system
\begin{equation}\label{eq-local}
A^{j}\hat{\mathbf{u}}^{j}=\bm{f}^{j},
\end{equation}
where $A^j$, $\bm{f}^j$, and $\hat{\mathbf{u}}^j$ are respectively  the local system matrix, local load vector, and local unknown solution. All of those are complex-valued --- i.e.
\begin{equation*}
A^j:=A^{j,\re}+iA^{j,\im},\quad\hat{\mathbf{u}}^j:=\hat{\mathbf{u}}^{j,\re}+i\hat{\mathbf{u}}^{j,\im}, \quad\bm{f}^j:=\bm{f}^{j,\re}+i\bm{f}^{j,\im}.
\end{equation*}

Similar to the real-valued case, the local system \eqref{eq-local} is singular due to the lack of boundary conditions. To address this, we decompose $\hat{\mathbf{u}}^{j}$ into two parts --- viz. the interface part $\hat{\mathbf{u}}^{j}_{\Gamma}$, where  nodes are shared by two or more adjacent subdomains, and the interior part $\hat{\mathbf{u}}^{j}_{I}$ belonging to the subdomain $D_j$.
As a result, the Eq.~\eqref{eq-local} can be written as
\begin{equation*}
  \begin{bmatrix}A_{II}^{j,\re}+iA_{II}^{j,\im} & A_{I\Gamma }^{j,\re}+iA_{I\Gamma }^{j,\im}\\
\\
A_{\Gamma I}^{j,\re}+iA_{\Gamma I}^{j,\im}& A_{\Gamma \Gamma }^{j,\re}+iA_{\Gamma \Gamma }^{j,\im} \\
\end{bmatrix}
\begin{Bmatrix}\hat{\mathbf{u}}^{j,\re}_{I}+i\hat{\mathbf{u}}^{j,\im}_{I}
\\
\\ \hat{\mathbf{u}}^{j,\re}_{\Gamma} +i\hat{\mathbf{u}}^{j,\im}_{\Gamma}\end{Bmatrix}
=
\begin{Bmatrix} \bm{f}^{j,\re}_{I}+i\bm{f}^{j,\im}_{I} \\\\ \bm{f}^{j,\re}_{\Gamma}+i\bm{f}^{j,\im}_{\Gamma}\end{Bmatrix}.
\end{equation*}
The system shows that once the interface unknowns $\hat{\mathbf{u}}^{j}_{\Gamma}$ are determined, the interior unknowns $\hat{\mathbf{u}}^{j}_{I}$ can be determined by solving  the interior problem
\begin{equation*}
    \left[A_{II}^{j,\re}+iA_{II}^{j,\im}\right] \left\{ \hat{\mathbf{u}}^{j,\re}_{I}+i\hat{\mathbf{u}}^{j,\im}_{I}\right\}=
 \left\{ \bm{f}^{j,\re}_{I}+i\bm{f}^{j,\im}_{I}\right\} - \left[A_{I\Gamma }^{j,\re}+iA_{I\Gamma }^{j,\im}\right]
\left\{\hat{\mathbf{u}}^{j,\re}_{\Gamma} +i\hat{\mathbf{u}}^{j,\im}_{\Gamma}\right\}.
\end{equation*}

Let $S_j$ and $F_j$ be respectively the local Schur complement matrix and the corresponding right-hand side vector --- i.e.
\begin{align*}
S_j&= \left[A_{\Gamma \Gamma }^{j,\re}+iA_{\Gamma \Gamma }^{j,\im}\right] - \left[A_{\Gamma I}^{j,\re}+iA_{\Gamma I}^{j,\im}\right] \left[A_{II}^{j,\re}+iA_{II}^{j,\im}\right]^{-1} \left[A_{I\Gamma }^{j,\re}+iA_{I\Gamma }^{j,\im}\right],\\
F_j&=\left\{\bm{f}^{j,\re}_{\Gamma}+i\bm{f}^{j,\im}_{\Gamma}\right\}-\left[A_{\Gamma I}^{j,\re}+iA_{\Gamma I}^{j,\im}\right] \left[A_{II}^{j,\re}+iA_{II}^{j,\im}\right]^{-1} \left\{ \bm{f}^{j,\re}_{I}+i\bm{f}^{j,\im}_{I}\right\}.
\end{align*}

Then the global interface unknowns $\hat{\mathbf{u}}_{\Gamma}$ can be obtained by solving the global interface problem
\begin{equation}\label{interface-d}
S\hat{\mathbf{u}}_{\Gamma}=F,
\end{equation}
where
\begin{equation*}
    S=\sum_{j=1}^{N_s}R_j^{T}S_jR_j, \ ~ F=\sum_{j=1}^{N_s}R_j^{T}F_j.
\end{equation*}
The restriction matrix $R_j$ acts as a scatter operator that maps the global interface unknowns
$\hat{\mathbf{u}}_{\Gamma}$ to the local interface unknowns $\hat{\mathbf{u}}^{j}_{\Gamma}$, satisfying the relation $R_j\hat{\mathbf{u}}_{\Gamma}=\hat{\mathbf{u}}^{j}_{\Gamma}.$
This operator enables the consistent transfer of interface information from the global domain to each local subdomain within the domain decomposition framework.

Note that the primary computational cost of the construction of $S$ and $F$ arises from inverting the complex-valued matrices $A_{II}^{j,\re}+iA_{II}^{j,\im}$, $ j = 1,\ldots,N_s$.
In this work, we  investigate the application of the domain decomposition method to parametric dynamical systems.
Using the Fourier transform, these systems can be converted into complex-valued parametric elliptic problems.
For each random input $\bm\mu\in\tilde\Omega$, all the matrices mentioned are $\bm\mu$-dependent, which leads to significant computational challenges when solving the interface system.

\section{Domain Decomposition Based on Variable-Separation}
\label{sec-sddm}\setcounter{equation}{0}
In this section, we will first present an offline-online method for the interface problem of the complex-valued elliptic problem. This builds a relation between random input and global interface solution --- cf. Sections \ref{sec-sf}-\ref{sec-reduce model}.
After that we present separate representations of the subproblems solutions.
According to the Eq.~\eqref{interface-d}, the interface problem for the complex-valued elliptic equation can be written as 
\begin{equation*}
  S(\bm{\mu})\hat{\mathbf{u}}_{\Gamma}(\bm{\mu})=F(\bm{\mu}),
\end{equation*}
where
\begin{equation*}
S(\bm{\mu})=\sum_{j=1}^{N_s}R_j^{T}S_j(\bm{\mu})R_j,\ ~ F(\bm{\mu})=\sum_{j=1}^{N_s}R_j^{T}F_j(\bm{\mu}),
\end{equation*}
and
\begin{align}
S_j(\bm{\mu})&= \left[A_{\Gamma \Gamma }^{j,\re}(\bm{\mu})+iA_{\Gamma \Gamma }^{j,\im}(\bm{\mu})\right]\nn\\
& - \left[A_{\Gamma I}^{j,\re}(\bm{\mu})+iA_{\Gamma I}^{j,\im}(\bm{\mu})\right] \left[A_{II}^{j,\re}(\bm{\mu})+iA_{II}^{j,\im}(\bm{\mu})\right]^{-1} \left[A_{I\Gamma }^{j,\re}(\bm{\mu})+iA_{I\Gamma }^{j,\im}(\bm{\mu})\right],\label{eq-S-F-1}\\
F_j(\bm{\mu})&=\left\{\bm{f}^{j,\re}_{\Gamma}(\bm{\mu})+i\bm{f}^{j,\im}_{\Gamma}(\bm{\mu})\right\}\nn\\
&-\left[A_{\Gamma I}^{j,\re}(\bm{\mu})+iA_{\Gamma I}^{j,\im}(\bm{\mu})\right] \left[A_{II}^{j,\re}(\bm{\mu})+iA_{II}^{j,\im}(\bm{\mu})\right]^{-1} \left\{ \bm{f}^{j,\re}_{I}(\bm{\mu})+i\bm{f}^{j,\im}_{I}(\bm{\mu})\right\}.\label{eq-S-F-2}
\end{align}

To simplify the presentation of the method, we split the domain $D$ into two subdomains $D_1$ and $D_2$ and denote their interface by $\Gamma$ instead of $\Gamma_{12}$. Let us point out that two-subdomain analysis extends naturally to multiple number of subdomains. The corresponding interface problem for the complex-valued elliptic equation takes the following form
\begin{equation}\label{eq-spde-interface}
(S_1(\bm{\mu})+S_2(\bm{\mu}))\hat{\mathbf{u}}_{\Gamma}(\bm{\mu})=F_1(\bm{\mu})+F_2(\bm{\mu}).
\end{equation}

\subsection{Affine representations of $S$ and $F$}
\label{sec-sf}
To enhance the computational efficiency, we develop an affine reformulation strategy for \eqref{eq-spde-interface} that guarantees the affine parameter dependence of both $S(\bm{\mu})$ and $F(\bm{\mu})$. The implementation details of this approach are presented in this section.
Using the affine decomposition \eqref{eq-fourier-affine}, we write the parametric matrices in equations \eqref{eq-S-F-1}-\eqref{eq-S-F-2} as
\begin{equation}\label{eq-affmatrix}
\begin{aligned}
    A_{II}^{j,\re}(\bm{\mu})&=\sum_{n=1}^{m_{a_j}}\alpha^{jn}(\bm{\mu})A_{II}^{jn},~~ A_{II}^{j,\im}(\bm{\mu})=\gamma^{j}(\bm{\mu})M_{II}^{j},\\
    A_{I\Gamma}^{j,\re}(\bm{\mu})&=\sum_{n=1}^{m_{a_j}}\alpha^{jn}(\bm{\mu})A_{I\Gamma}^{jn},~~A_{I\Gamma}^{j,\im}(\bm{\mu})=\gamma^{j}(\bm{\mu})M_{I\Gamma}^{j},\\
    A_{\Gamma I}^{j,\re}(\bm{\mu})&=\sum_{n=1}^{m_{a_j}}\alpha^{jn}(\bm{\mu})A_{\Gamma I}^{jn},~~
    A_{\Gamma I}^{j,\im}(\bm{\mu})=\gamma^{j}(\bm{\mu})M_{\Gamma I}^{j},\\
    A_{\Gamma\Gamma}^{j,\re}(\bm{\mu})&=\sum_{n=1}^{m_{a_j}}\alpha^{jn}(\bm{\mu})A_{\Gamma\Gamma}^{jn},~~A_{\Gamma\Gamma}^{j,\im}(\bm{\mu})=\gamma^{j}(\bm{\mu})M_{\Gamma\Gamma}^{j},\\
\end{aligned}
\end{equation}
where the matrices $A_{II}^{jn}$, $A_{I\Gamma}^{jn}$, $A_{\Gamma I}^{jn}$, $A_{\Gamma\Gamma}^{jn}$ and $M_{II}^{j}$, $M_{I\Gamma}^{j}$, $M_{\Gamma I}^{j}$, $M_{\Gamma\Gamma}^{j}$ do not depend on $\bm{\mu}$.
Note that the first term of $S_j(\bm{\mu})$ has the structure affine with respect to $\bm{\mu}$ as required. Next, we  discuss how to obtain the affine approximation of the term 
\begin{equation*}
-\left[A_{\Gamma I}^{j,\re}(\bm{\mu})+iA_{\Gamma I}^{j,\im}(\bm{\mu})\right] \left[A_{II}^{j,\re}(\bm{\mu})+iA_{II}^{j,\im}(\bm{\mu})\right]^{-1} \left[A_{I\Gamma }^{j,\re}(\bm{\mu})+iA_{I\Gamma }^{j,\im}(\bm{\mu})\right]
\end{equation*}
in \eqref{eq-S-F-1}.

\begin{Steps}

\item Construct the low-rank approximation of
\begin{equation}\label{eq-X-def}
{X}^j(\bm{\mu}) = \left[A_{II}^{j,\re}(\bm{\mu})+iA_{II}^{j,\im}(\bm{\mu})\right]^{-1} \left[A_{I\Gamma }^{j,\re}(\bm{\mu})+iA_{I\Gamma }^{j,\im}(\bm{\mu})\right]
\end{equation}
such as
\begin{equation}\label{eq-X-appr}
{X}^j_{N}(\bm{\mu}):=\sum_{n=1}^{N_{S_j}}\left( \phi_n^{j,\re}(\bm{\mu}){X}_n^{j,\re} + i \phi_n^{j,\im}(\bm{\mu}){X}_n^{j,\im} \right).
\end{equation}

Indeed, under the assumption \eqref{eq-affmatrix}, one can rewrite \eqref{eq-X-def} as
\begin{equation*}
    \left(\sum_{n=1}^{m_{a_j}}\alpha^{jn}(\bm{\mu})A_{II}^{jn}+i\gamma^{j}(\bm{\mu})M_{II}^{j}\right)
X^j(\bm{\mu})
=
\left(\sum_{n=1}^{m_{a_j}}\alpha^{jn}(\bm{\mu})A_{I\Gamma}^{jn}+i\gamma^{j}(\bm{\mu})M_{I\Gamma}^{j}\right).
\end{equation*}

Let $n_{\Gamma}$ be the number of the interface unknowns. The unknown parametric matrix $X^j(\bm{\mu})$ can be expressed as $X^j(\bm{\mu})=[x_1^{j}(\bm{\mu}), x_2^{j}(\bm{\mu}), \ldots,x_{n_{\Gamma}}^{j}(\bm{\mu})]$,
and each element $x_k^{j}(\bm{\mu}), k=1, \ldots, n_\Gamma$ can be obtained by solving
\begin{equation}\label{eq-stepass}
\left(\sum_{n=1}^{m_{a_j}}\alpha^{jn}(\bm{\mu})A_{II}^{jn}+i\gamma^{j}(\bm{\mu})M_{II}^{j}\right)x_k^{j}(\bm{\mu})
=\delta_k^{j}(\bm{\mu}),
\end{equation}
where $\delta_k^{j}(\bm{\mu})$ denotes the $k$-th column of the matrix $\left(\sum_{n=1}^{m_{a_j}}\alpha^{jn}(\bm{\mu})A_{I\Gamma}^{jn}+i\gamma^{j}(\bm{\mu})M_{I\Gamma}^{j}\right)$.

Applying Algorithm \ref{algorithm-vs-sSAS}, we can  derive the low-rank approximation for each component $x_k^{j}(\bm{\mu})$.
Then, the approximate solution $X_N^j(\bm{\mu})$ in \eqref{eq-X-appr} is constructed by rearranging the low-rank approximations of $\{x_k^{j}(\bm{\mu})\}_{k=1}^{n_\Gamma}$.

\item Assemble affine expression for
\begin{equation*}		
\mathcal{X}^j(\bm{\mu})=
-\left[A_{\Gamma I}^{j,\re}(\bm{\mu})+iA_{\Gamma I}^{j,\im}(\bm{\mu})\right] \left[A_{II}^{j,\re}(\bm{\mu})+iA_{II}^{j,\im}(\bm{\mu})\right]^{-1} \left[A_{I\Gamma }^{j,\re}(\bm{\mu})+iA_{I\Gamma }^{j,\im}(\bm{\mu})\right].
\end{equation*}

Using the low-rank representation of $X^j(\bm{\mu})$ and the Eq.~\eqref{eq-affmatrix}, we write
\begin{align}
\mathcal{X}^j(\bm{\mu})&\approx -\Big[\sum_{n=1}^{m_{a_j}}\alpha^{jn}(\bm{\mu})A_{\Gamma I}^{jn}+i\gamma^{j}(\bm{\mu})M_{\Gamma I}^{j}\Big] \Big[\sum_{n=1}^{N_{S_j}}\left( \phi_n^{j,\re}(\bm{\mu}){X}_n^{j,\re} + i \phi_n^{j,\im}(\bm{\mu}){X}_n^{j,\im} \right)\Big]\nn\\
&=\sum_{n=1}^{(m_{a_j}+1)N_{S_j}}\left( \eta_n^{j,\re}(\bm{\mu})\mathcal{X}_n^{j,\re} + i\eta_n^{j,\im}(\bm{\mu})\mathcal{X}_n^{j,\im} \right),\label{eq-x-def}
\end{align}
where
\begin{align*}
  \eta_n^{j,\re}(\bm{\mu})&=\alpha^{jn_1}(\bm{\mu})\phi_{n_2}^{j,\re}(\bm{\mu}),\
\mathcal{X}_n^{j,\re}=-A_{\Gamma I}^{jn_1}{X}_{n_2}^{j,\re},\\
\eta_n^{j,\im}(\bm{\mu})&=\alpha^{jn_1}(\bm{\mu})\phi_{n_2}^{j,\im}(\bm{\mu}),
\ \mathcal{X}_n^{j,\im}=-A_{\Gamma I}^{jn_1}{X}_{n_2}^{j,\im},\quad 
n=1,2,\ldots,m_{a_j}N_{S_j},  
\end{align*}
and
\begin{align*}
\eta_n^{j,\re}(\bm{\mu})&=\gamma^{j}(\bm{\mu})\phi_{n_2}^{j,\im}(\bm{\mu}),\
\mathcal{X}_n^{j,\re}=M_{\Gamma I}^{j} {X}_{n_2}^{j,\im}, \\
\eta_n^{j,\im}(\bm{\mu})&=\gamma^{j}(\bm{\mu})\phi_{n_2}^{j,\re}(\bm{\mu}),
\mathcal{X}_n^{j,\im}=-M_{\Gamma I}^{j}{X}_{n_2}^{j,\re},\quad 
n=m_{a_j}N_{S_j}+1,\ldots,(m_{a_j}+1)N_{S_j}, 
\end{align*}
where $n_1=1,2,\ldots,m_{a_j}$ and $n_2=1,2,\ldots,N_{S_j}$.
\end{Steps}

Subsequently, the low-rank representation of $S_j(\bm{\mu})$ takes the form
\begin{align*}
S_j(\bm{\mu})&=
 \left[A_{\Gamma \Gamma }^{j,\re}(\bm{\mu})+iA_{\Gamma \Gamma }^{j,\im}(\bm{\mu})\right]
+ \mathcal{X}^j(\bm{\mu})
\\&\approx
\nonumber
\sum_{n=1}^{m_{a_j}}\alpha^{jn}(\bm{\mu})A_{\Gamma\Gamma}^{jn}+i\gamma^{j}(\bm{\mu})M_{\Gamma\Gamma}^{j}
+ \sum_{n=1}^{(m_{a_j}+1)N_{S_j}}\left( \eta_n^{j,\re}(\bm{\mu})\mathcal{X}_n^{j,\re} + i\eta_n^{j,\im}(\bm{\mu})\mathcal{X}_n^{j,\im}\right),
\end{align*}
so that the low-rank representation of $S(\bm{\mu})$ is
\begin{equation}\label{eq-S-repre}
S(\bm{\mu})= S_1(\bm{\mu})+S_2(\bm{\mu})
\approx
\sum_{n=1}^{m_S}\left( \hat{\eta}_n^{\re}(\bm{\mu})\mathcal{\hat{X}}_n^{\re} + i\hat{\eta}_n^{\im}(\bm{\mu})\mathcal{\hat{X}}_n^{\im}\right),
\end{equation}
where the sum term comes from stacking the variables and sorting the corresponding indices via a single one, and  $m_S=m_{a_1} +m_{a_2} +  (m_{a_1}+1)N_{S_1}+ (m_{a_2}+1)N_{S_2}$.
Note that $m_S$ is the number of the real part of $S(\bm{\mu})$, and the number of its imaginary part is equal to $2 + (m_{a_1} + 1)N_{S_1} + (m_{a_2} + 1)N_{S_2}$, which is less than or equal to $m_S$. For convenience, we also use $m_S$ to denote the number of the imaginary part, since for the imaginary part, the terms with the index
$n=3 + (m_{a_1} + 1)N_{S_1} + (m_{a_2} + 1)N_{S_2},\ldots,m_S$ can be set to $0$.
Algorithm~\ref{algorithm-S} outlines the assembling process of $S(\bm{\mu})$.~%
\begin{algorithm}[!tb]
\caption{The assemble process of $S(\bm{\mu})$}
\label{algorithm-S}
\begin{algorithmic}
\REQUIRE{The complex-valued matrices $A_{II}^{j}(\bm{\mu})$, $A_{I\Gamma}^{j}(\bm{\mu})$, $A_{\Gamma I}^{j}(\bm{\mu})$, $A_{\Gamma\Gamma}^{j}(\bm{\mu})$, $j=1,2$.}  \\
\ENSURE{The low-rank representation $S(\bm{\mu})=\sum_{n=1}^{m_S}\left( \hat{\eta}_n^{\re}(\bm{\mu})\mathcal{\hat{X}}_n^{\re} + i\hat{\eta}_n^{\im}(\bm{\mu})\mathcal{\hat{X}}_n^{\im}\right)$.}
\begin{enumerate}
    \item
Get the approximation of
\begin{equation*}
  {X}^j(\bm{\mu}) = \left[A_{II}^{j,\re}(\bm{\mu})+iA_{II}^{j,\im}(\bm{\mu})\right]^{-1} \left[A_{I\Gamma }^{j,\re}(\bm{\mu})+iA_{I\Gamma }^{j,\im}(\bm{\mu})\right]  
\end{equation*}
by solving Eq.~\eqref{eq-stepass} with Algorithm \ref{algorithm-vs-sSAS};
\item Assemble the affine expression of 
\begin{equation*}
  -\left[A_{\Gamma I}^{j,\re}(\bm{\mu})+iA_{\Gamma I}^{j,\im}(\bm{\mu})\right] \left[A_{II}^{j,\re}(\bm{\mu})+iA_{II}^{j,\im}(\bm{\mu})\right]^{-1} \left[A_{I\Gamma }^{j,\re}(\bm{\mu})+iA_{I\Gamma }^{j,\im}(\bm{\mu})\right]  
\end{equation*}
by \eqref{eq-x-def};
\item
Assemble $S_j(\bm{\mu})$ based on \eqref{eq-affmatrix} and the expression derived in Step 2;
\item
Assemble the low-rank representation of $S(\bm{\mu})$ by \eqref{eq-S-repre}.
\end{enumerate}
\end{algorithmic}
\end{algorithm}

Similarly, with the assumption \eqref{eq-fourier-affine} of affine decomposition, we have
\begin{equation}\label{eq-affvector}
\begin{aligned}
\bm{f}_I^{j,\re}(\bm{\mu})=\sum_{n=1}^{m_{b_j}}\beta^{jn,\re}(\bm{\mu})\bm{f}_I^{jn,\re}, ~~
\bm{f}_I^{j,\im}(\bm{\mu})=\sum_{n=1}^{m_{b_j}}\beta^{jn,\im}(\bm{\mu})\bm{f}_I^{jn,\im},\\
\bm{f}_{\Gamma}^{j,\re}(\bm{\mu})=\sum_{n=1}^{m_{b_j}}\beta^{jn,\re}(\bm{\mu})\bm{f}_{\Gamma}^{jn,\re},~~
\bm{f}_{\Gamma}^{j,\im}(\bm{\mu})=\sum_{n=1}^{m_{b_j}}\beta^{jn,\im}(\bm{\mu})\bm{f}_{\Gamma}^{jn,\im},
\end{aligned}
\end{equation}
where $\bm{f}_I^{jn,\re}$, $\bm{f}_I^{jn,\im}$, $\bm{f}_{\Gamma}^{jn,\re}$ and $\bm{f}_{\Gamma}^{jn,\im}$ are independent of $\bm{\mu}$.
Following the assemble strategy of $S(\bm{\mu})$, one can derive the low-rank representation for $F(\bm{\mu})$ as follows
\begin{equation}\label{eq-F-repre}
F(\bm{\mu})= F_1(\bm{\mu})+F_2(\bm{\mu})
\approx\sum_{n=1}^{m_F}\left( \hat{\rho}_n^{\re}(\bm{\mu}){\hat{F}}_n^{\re} + i\hat{\rho}_n^{\im}(\bm{\mu}){\hat{F}}_n^{\im}\right),
\end{equation}
where
\begin{equation*}
m_F = m_{b_1}+m_{b_2}+(m_{a_1}+1)N_{F_1}+(m_{a_2}+1)N_{F_2}.
\end{equation*}
The detailed procedure of assembling $F(\bm{\mu})$ is presented in Algorithm \ref{algorithm-F}.

\begin{algorithm}
\caption{The assemble process of $F(\bm{\mu})$ }
\label{algorithm-F}
\begin{algorithmic}
\REQUIRE{The complex-valued matrices $A_{II}^{j}(\bm{\mu})$, $A_{\Gamma I}^{j}(\bm{\mu})$ and vectors $\bm{f}_I^{j}(\bm{\mu})$, $\bm{f}_{\Gamma}^{j}(\bm{\mu})$, $j=1,2$.}
\ENSURE{The low-rank representation $F(\bm{\mu})=\sum_{n=1}^{m_F}\left( \hat{\rho}_n^{\re}(\bm{\mu}){\hat{F}}_n^{\re} + i\hat{\rho}_n^{\im}(\bm{\mu}){\hat{F}}_n^{\im}\right)$.}
\begin{enumerate}
    \item
Get the approximation of $\left[A_{II}^{j,\re}(\bm{\mu})+iA_{II}^{j,\im}(\bm{\mu})\right]^{-1} \left\{\bm{f}_{I}^{j,\re}(\bm{\mu}) + i \bm{f}_{I}^{j,\im}(\bm{\mu})\right\}$ by Algorithm \ref{algorithm-vs-sSAS};
\item
Assemble the affine expression of 
\begin{equation*}
    -\left[A_{\Gamma I}^{j,\re}(\bm{\mu})+iA_{\Gamma I}^{j,\im}(\bm{\mu})\right] \left[A_{II}^{j,\re}(\bm{\mu})+iA_{II}^{j,\im}(\bm{\mu})\right]^{-1} \left\{\bm{f}_{I}^{j,\re}(\bm{\mu}) + i \bm{f}_{I}^{j,\im}(\bm{\mu})\right\};
\end{equation*}
\item
Assemble $F_j(\bm{\mu})$ based on \eqref{eq-affvector} and the expression derived in Step 2;
\item
Assemble the low-rank representation of $F(\bm{\mu})$ by \eqref{eq-F-repre}.
\end{enumerate}
\end{algorithmic}
\end{algorithm}

\subsection{Reduced model for the interface problem}
\label{sec-reduce model}
Based on the separate approximations of $S(\bm{\mu})$ and $F(\bm{\mu})$ derived in subsection \ref{sec-sf}, the interface problem \eqref{eq-spde-interface} can be expressed as
\begin{equation}\label{eq-Intaff}
\sum_{n=1}^{m_S}\left(\hat{\eta}_n^{\re}(\bm{\mu})\mathcal{\hat{X}}_n^{\re} + i\hat{\eta}_n^{\im}(\bm{\mu})\mathcal{\hat{X}}_n^{\im}\right)\hat{\mathbf{u}}_{\Gamma}(\bm{\mu})=\sum_{n=1}^{m_F}\left( \hat{\rho}_n^{\re}(\bm{\mu}){\hat{F}}_n^{\re} + i\hat{\rho}_n^{\im}(\bm{\mu}){\hat{F}}_n^{\im}\right),
\end{equation}
which is a linear complex-valued system with $n_\Gamma$ unknowns.
Although solving \eqref{eq-Intaff} requires much less computation effort than that of the original model \eqref{eq-spde-interface}, this may not be a very small-scale problem.
This persistence stems from the fact that the reformulated model \eqref{eq-Intaff} inherits the discretization-dependent characteristics of the original full-order system.
Therefore, we again employ the VS method from Section \ref{sec-vs} to derive the following reduced model representation of the interface problem:
\begin{equation}\label{eq-interu}
\hat{\mathbf{u}}_{\Gamma}(\bm{\mu})\approx\sum_{k=1}^{N_{\Gamma}}\left(\zeta_{\Gamma}^{k,\re}(\bm{\mu}) \mathbf{c}_{\Gamma}^{k,\re}+i\zeta_{\Gamma}^{k,\im}(\bm{\mu}) \mathbf{c}_{\Gamma}^{k,\im}\right),
\end{equation}
where $N_{\Gamma}$ is the number of the separate terms, $\bm{c}_{\Gamma}^{k,\re}$ and $\bm{c}_{\Gamma}^{k,\im}$, $k=1,\ldots,N_{\Gamma}$ are $n_\Gamma$-dimension vectors.
Using this separate representation for new random inputs, we can directly obtain the solution to the interface problem without solving the Eq. \eqref{eq-spde-interface} or the Eq. \eqref{eq-Intaff} again.
This approach has the potential of very high efficiency, since it only requires computing the coefficients $\zeta_{\Gamma}^{k,\re}(\bm{\mu})$ and $\zeta_{\Gamma}^{k,\im}(\bm{\mu})$ for new parameter values while reusing the precomputed spatial basis $\mathbf{c}_{\Gamma}^{k,\re}$ and $\mathbf{c}_{\Gamma}^{k,\im}$.

\subsection{Reduced model for subdomain problems}
\label{sec-subproblems}
As discussed in  Section \ref{sec-ddm}, once having the solution $ \hat{\mathbf{u}}_{\Gamma}^j(\bm{\mu})$ of the interface problem, the interior solution $ \hat{\mathbf{u}}_{I}^j(\bm{\mu})$ in the subdomain $D_j$ can be obtained by solving the complex-valued equation
\begin{align*}
&\left[A_{II}^{j,\re}(\bm{\mu})+iA_{II}^{j,\im}(\bm{\mu})\right]
\left\{ \hat{\mathbf{u}}^{j,\re}_{I}(\bm{\mu})+i\hat{\mathbf{u}}^{j,\im}_{I}(\bm{\mu})\right\}\\
&\quad\quad=\left\{ \bm{f}^{j,\re}_{I}(\bm{\mu})+i\bm{f}^{j,\im}_{I}(\bm{\mu})\right\} - \left[A_{I\Gamma }^{j,\re}(\bm{\mu})+iA_{I\Gamma }^{j,\im}(\bm{\mu})\right]
\left\{\hat{\mathbf{u}}^{j,\re}_{\Gamma}(\bm{\mu}) +i\hat{\mathbf{u}}^{j,\im}_{\Gamma}(\bm{\mu})\right\}.
\end{align*}
To reduce the computational cost, we apply the VS method in order to construct an efficient surrogate model.
Using the affine decomposition of the matrices $A_{II}^{j,\re}(\bm{\mu})$, $A_{II}^{j,\im}(\bm{\mu})$, $A_{IT}^{j,\re}(\bm{\mu})$, $A_{IT}^{j,\im}(\bm{\mu})$ in \eqref{eq-affmatrix} and the vectors $\bm{f}^{j,\re}_{I}(\bm{\mu})$, $\bm{f}^{j,\im}_{I}(\bm{\mu})$ in  \eqref{eq-affvector} and the affine surrogate model of the interface solutions $\hat{\mathbf{u}}^{j,\re}_{\Gamma}(\bm{\mu})$, $\hat{\mathbf{u}}^{j,\im}_{\Gamma}(\bm{\mu})$ in equation \eqref{eq-interu}, we get
\begin{equation*}
\hat{\mathbf{u}}_{I}^j(\bm{\mu})\approx\sum_{k=1}^{N_{I}}\left(\zeta_{I}^{jk,\re}(\bm{\mu}) \mathbf{c}_{I}^{jk,\re}+i\zeta_{I}^{jk,\im}(\bm{\mu}) \mathbf{c}_{I}^{jk,\im}\right),
\end{equation*}
where $N_{I}$ is the number of the separate terms, $\bm{c}_{I}^{jk,\re}$ and $\bm{c}_{I}^{jk,\im}$, $k=1,\ldots,N_{I}$ are vectors independent of $\bm{\mu}$.

In summary, for time-dependent problems, after computing the corresponding solutions using the preconstructed surrogate models for both the interface problem and subproblems with new parameter inputs, we can obtain the temporal solution through an inverse Fourier transform. The entire computational process is independent of the spatial discretization of the original system and requires no time-stepping iterations, thus demonstrating remarkable efficiency.

Notably, the online stage of the proposed approach can be much more efficient, as the surrogate models of the interface problem and subproblems (as well as their real and imaginary parts) allow parallel computation across various frequency values.

\begin{remark}
The proposed method can be extended to non-affine, nonlinear, and multi-physical parameterized problems. For the non-affine case, the VS method \cite{Li2017a} for a multi-variable function and the empirical interpolation method can be used to obtain the affine form. In multi-physical parameterized problems, we can use the domain decomposition to reduce a global multi-physical problem to a local single-physical problem. Then we perform the VS method for each local problem. To obtain the Fourier transform of the nonlinear case, the convolution property will be used.
\end{remark}

\section{Numerical Experiments}
\label{sec-numerical examples}
Here we present three numerical examples. The corresponding numerical results demonstrate the applicability and efficiency of the method proposed.  In Section \ref{exa:1}, the method is applied to the heat equation. After that we apply it to the reaction diffusion equation in Sections \ref{exa:2} and \ref{exa:3} with different settings.
All numerical experiments  were run in Python on a Dell workstation with Intel Xeon W-2295 CPU @3.00GHz and 64GB of RAM.

To quantify the accuracy of the proposed method, the average relative error $\epsilon_u$ for the solution of the original dynamical systems is defined by
\begin{equation}\label{eq-meanerror}
\epsilon_{u}:=\frac{1}{M}\sum_{k=1}^{M}\frac{||{u}(\bm{x},t;\bm{\xi}_k)-{u}_{N}(\bm{x},t;\bm{\xi}_k)||_{L^2([0,T];\mathcal{V}_h)}}{||{u}(\bm{x},t;\bm{\xi}_k)||_{L^2([0,T];\mathcal{V}_h)}},
\end{equation}
where $M$ is the number of samples, ${u}_{N}(\bm{x},t;\bm{\xi}_k)$ the approximate solution obtained by the proposed method which is denoted FT-DD-VS, and ${u}(\bm{x},t;\bm{\xi}_k)$ the reference solution.
In the following  examples, the reference solution is calculated by FEM in space and the backward Euler scheme in time, denoted as FEM-BE.

Furthermore, we define the average relative error $\epsilon_{\hat{\mathbf{u}}}$ for the complex-valued equations as
\begin{equation}\label{eq-meanerror-2}
\epsilon_{\hat{\mathbf{u}}}:=\frac{1}{M}\sum_{k=1}^{M}\frac{|\hat{\mathbf{u}}(\bm{\mu}_k)-\hat{\mathbf{u}}_{N}(\bm{\mu}_k)|}{|\hat{\mathbf{u}}(\bm{\mu}_k)|},
\end{equation}
where $\hat{\mathbf{u}}_{N}(\bm{\mu}_k)$ is the approximate solution obtained by the VS method and $\hat{\mathbf{u}}(\bm{\mu}_k)$ the reference solution determined by FEM.

\subsection{Heat equation}
\label{exa:1}
On the domain $D=[0,1]\times [0,1]$, we consider the following heat equation:
\begin{equation*}
    \frac{\partial u}{\partial t}(\bm{x},t;\bm{\xi})-\nabla\cdot (c(\bm{x};\bm{\xi}) \nabla u(\bm{x},t;\bm{\xi}))=f(\bm{x},t;\bm{\xi}),
\end{equation*}
where
\begin{equation}\label{eq:exa1}
    \begin{aligned}
    u(\bm x,0;\bm\xi)&=0, \ (\bm x,\bm\xi)\in D\times\Omega,  \\
u(\bm x,t;\bm\xi)&=0,\ (\bm x,t;\bm\xi)\in \partial D\times[0,T]\times\Omega.
    \end{aligned}
\end{equation}
The analytical solution of this problem is
\begin{equation*}
    u(\bm{x},t;\bm{\xi}) = \frac{1}{\pi}\frac{t}{t^2+1}\sin{\pi {x}_1}\sin{\pi {x}_2},
\end{equation*}
so that the source function $f$ has the form
\begin{equation*}
(\bm{x},t;\bm{\xi}) = \left(\frac{1}{\pi}\frac{1-t^2}{(t^2+1)^2} + 2\pi \frac{t}{t^2+1}c(\bm{x};\bm{\xi})\right)\sin{\pi {x}_1}\sin{\pi {x}_2}.
\end{equation*}

The original domain $D$ is divided into the subdomains $D_1 = [0,0.5]\times [0,1]$ and $D_2 = [0.5,1]\times [0,1]$, and the coefficient $c(\bm{x};\bm{\xi})$ is defined as
\begin{equation*}
c(\bm{x};\bm{\xi}) = \left\{
\begin{aligned}
\xi_1,  & \ \bm{x}\in D_1,\\
2\xi_2,  &\ \bm{x}\in D_2.
\end{aligned}
\right. 
\end{equation*}
Therefore, the source function has different representations in subdomains $D_1$ and $D_2$.
Note that $\bm{\xi}:=({\xi}_1,{\xi}_2)\in[1,2]^2$ and $T=1$.
The reference solution is obtained by FEM in space with the mesh size $h=0.02$ and the backward Euler scheme in time with the step size $\tau=5\times 10^{-4}$. Besides, $|\Xi|$ denotes the cardinality of the sample set $\Xi$ and is set to $|\Xi|=10$ for the offline stage in Algorithm \ref{algorithm-vs-sSAS}.
In this problem, we take $\omega^*=20$ and $N_\omega = 20$, which suffices for determining an effective approximation.

\begin{figure}[!tb]
    \centering
    \includegraphics[width=7cm]{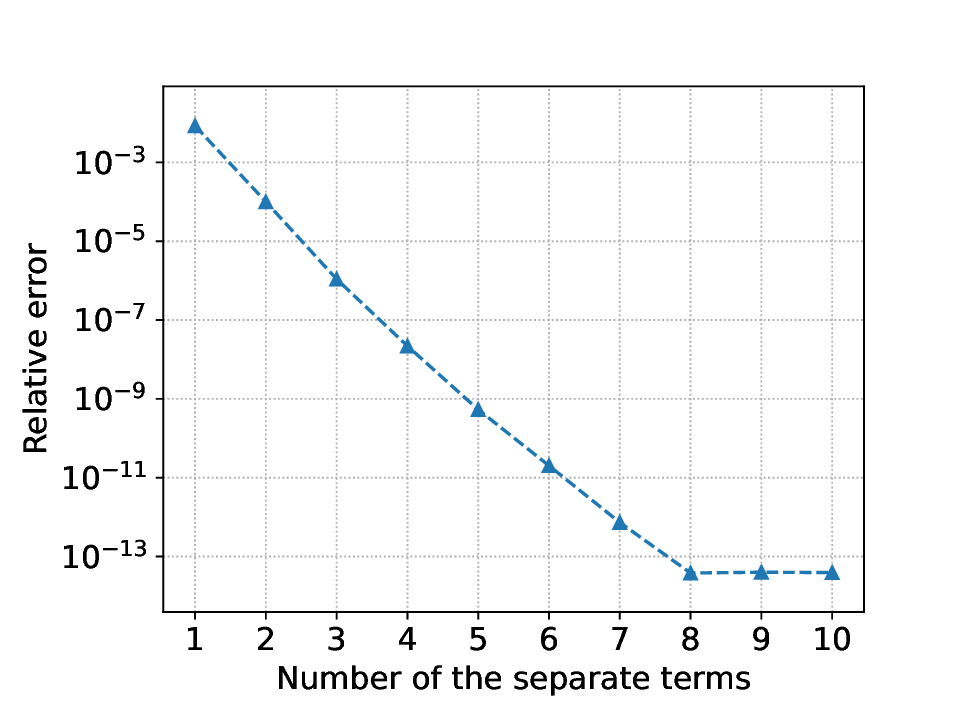 }
    \caption{Average relative error versus the number of separate terms $N_{S_1}$.}
    \label{fig-heat-11}
\end{figure}

First, we quantify the assembly errors for $S$ and $F$ in Step 1, which coincide with their affine approximation errors.
Taking $S_1$ as an example, Fig. \ref{fig-heat-11} plots the average relative error versus the number of separate terms $N_{S_1}$, computed as \eqref{eq-meanerror-2} with $10^3$ samples.
Note that the error diminishes monotonically as $N_{S_1}$ increases. This implies that the VS method effectively approximates complex-valued problems. To construct an efficient surrogate model for the interface problem, we select $N_{S_1}=4$, ensuring the relative approximation error of $S_1$ remains below $10^{-6}$.
Correspondingly, we set $N_{S_2} = 4$, $N_{F_1} = 4$, and $N_{F_2} = 4$.

\begin{figure}[!tbh]
\centering
\begin{minipage}{0.45\textwidth}
\centering
\includegraphics[width=6.5cm]{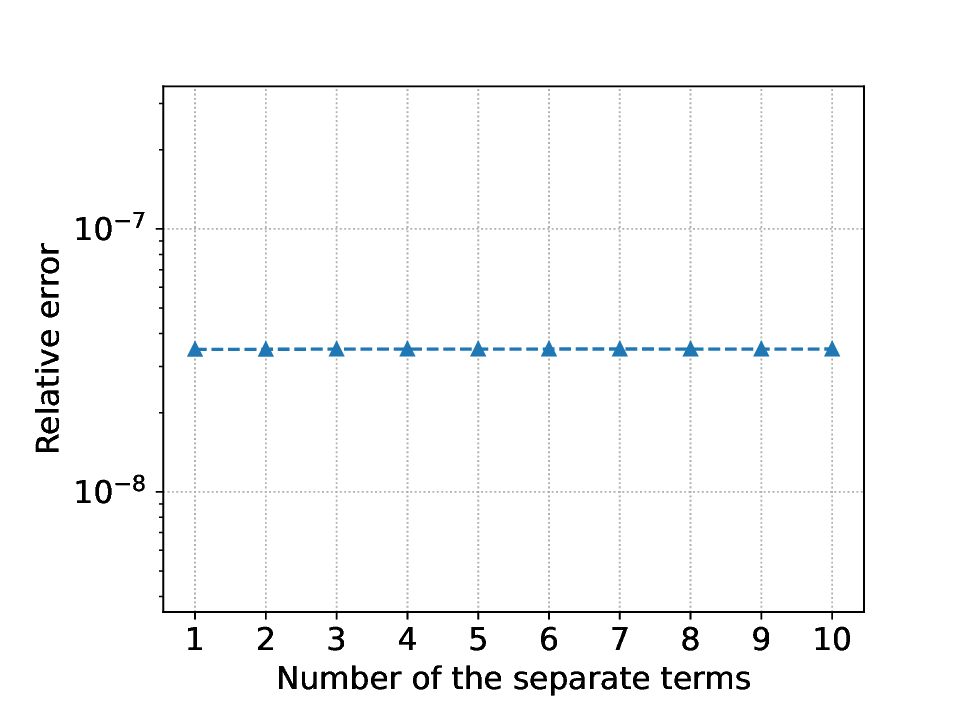 }\\
\scriptsize{(a) Interface problem}
\end{minipage}
\begin{minipage}{0.45\textwidth}
\centering
\includegraphics[width=6.5cm]{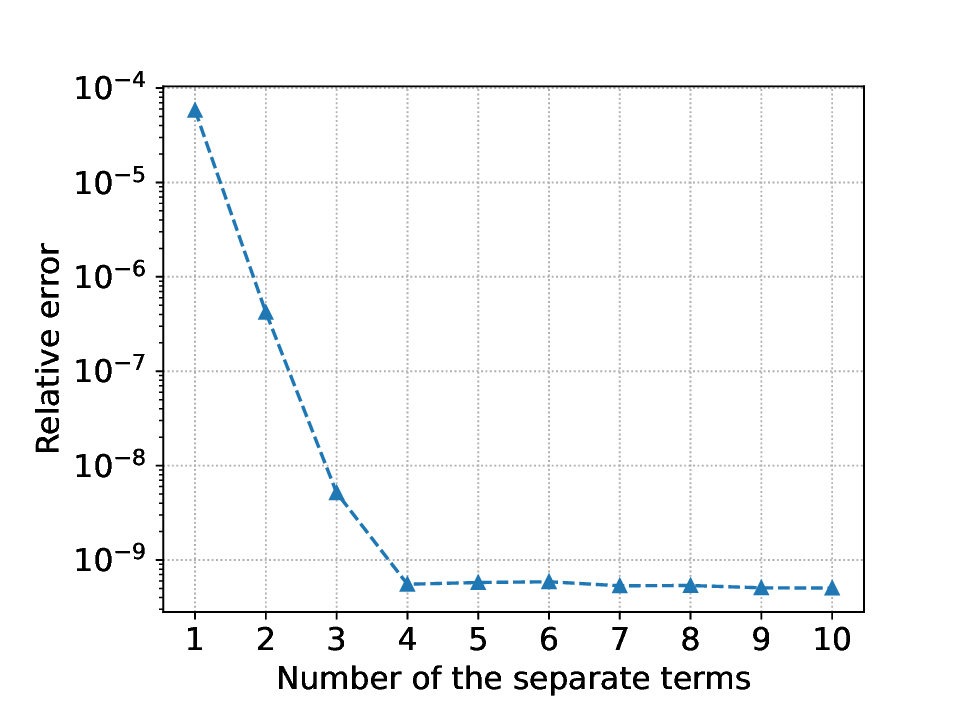 }\\
\scriptsize{(b) Subdomain problems}
\end{minipage}
\caption{Average relative error versus the number of separate terms.}
\label{fig-heat-1}
\end{figure}

Based on the affine approximations of $S$ and $F$, Fig. \ref{fig-heat-1} shows the average relative error versus the number of separate terms for interface problem $N_\Gamma$  and subdomain problems $N_I$.
The errors $\epsilon_{\mathbf{u}_{\Gamma}}$ and $\epsilon_{\mathbf{u}_I}$ are calculated by formula \eqref{eq-meanerror-2} with $10^3$ samples.

We find that the relative error of subdomain problems decreases rapidly and stabilizes when $N_I\geq4$. The results demonstrate that the VS method yields accurate approximations for both interface and subdomain solutions, a capability critical for ensuring the reliability of our proposed method in parametric dynamical systems.

\begin{figure}[!tbh]
\centering
\begin{minipage}{0.45\textwidth}
\centering
\includegraphics[width=6.5cm]{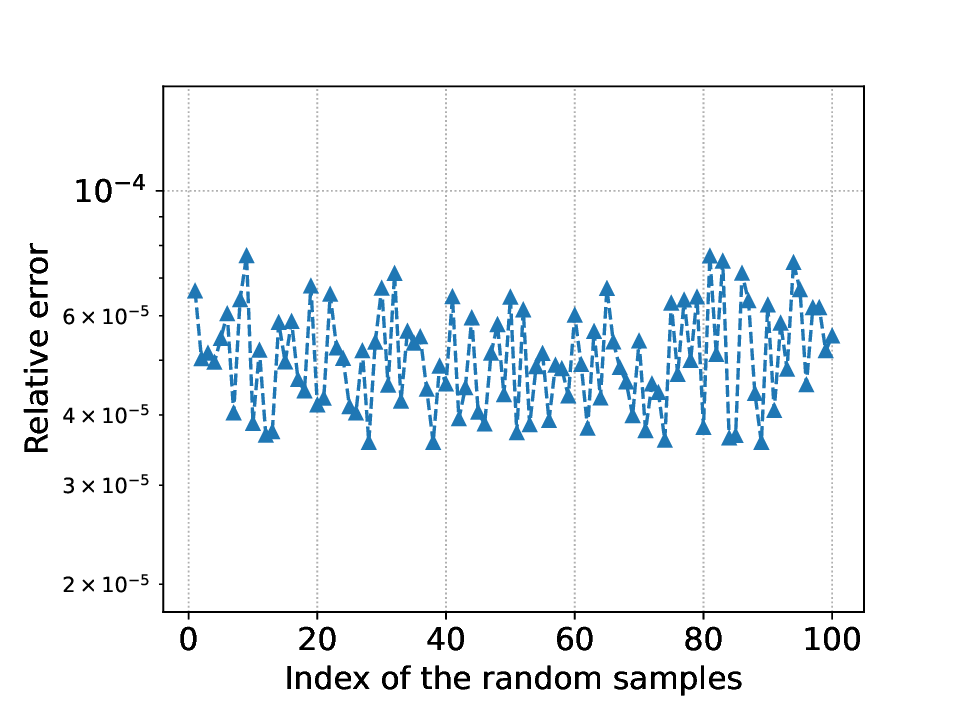 }\\
\scriptsize{(a) Errors for the first $100$ samples}
\end{minipage}
\begin{minipage}{0.45\textwidth}
\centering
\includegraphics[width=6.5cm]{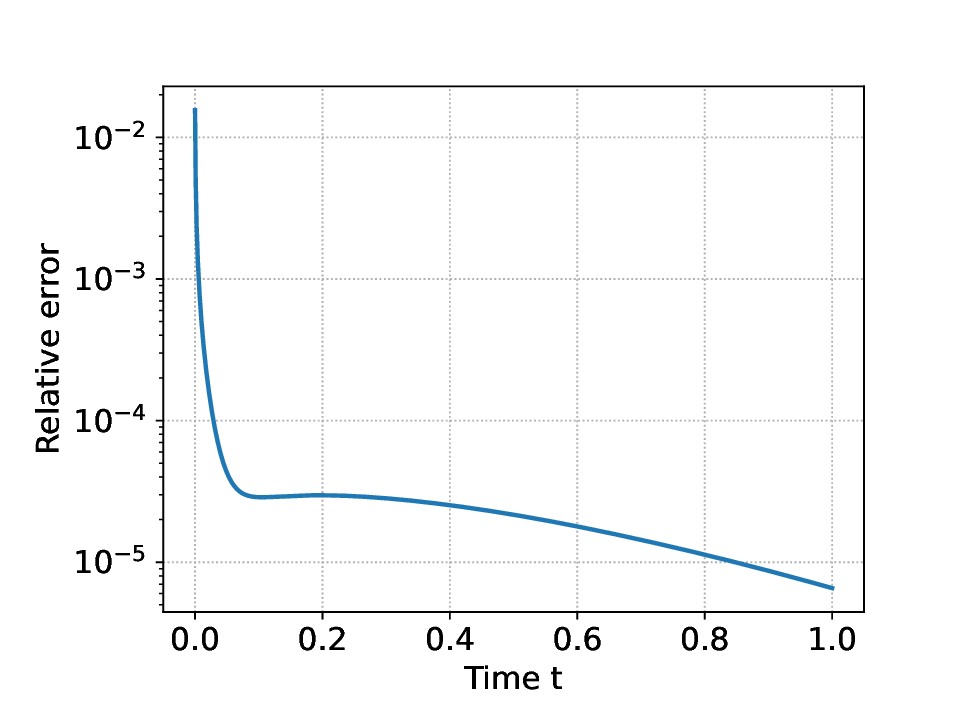 }\\
\scriptsize{(b) Average relative errors versus the time}
\end{minipage}
\caption{Comparison of the errors.}
\label{fig-heat-2}
\end{figure}

Fig. \ref{fig-heat-1}  shows that the number of separate terms $N_{\Gamma}=1$ for the interface problem and $N_I=4$ for the subdomain problems yields a surrogate model with the desired relative error while maintaining computational efficiency.
Under this configuration ($N_\Gamma=1$ and $N_I=4$), we compute the relative error of the original dynamical system by \eqref{eq-meanerror}. Fig. \ref{fig-heat-2} (a) illustrates the relative errors of the first 100 samples. Along with the temporal evolution, the relative error decreases --- cf. Fig. \ref{fig-heat-2} (b). Notice that the proposed method can provide a good approximation for the solution of the original problem.

\begin{figure}[!tbh]
\centering
\begin{minipage}{0.45\textwidth}
\centering
\includegraphics[width=6.5cm]{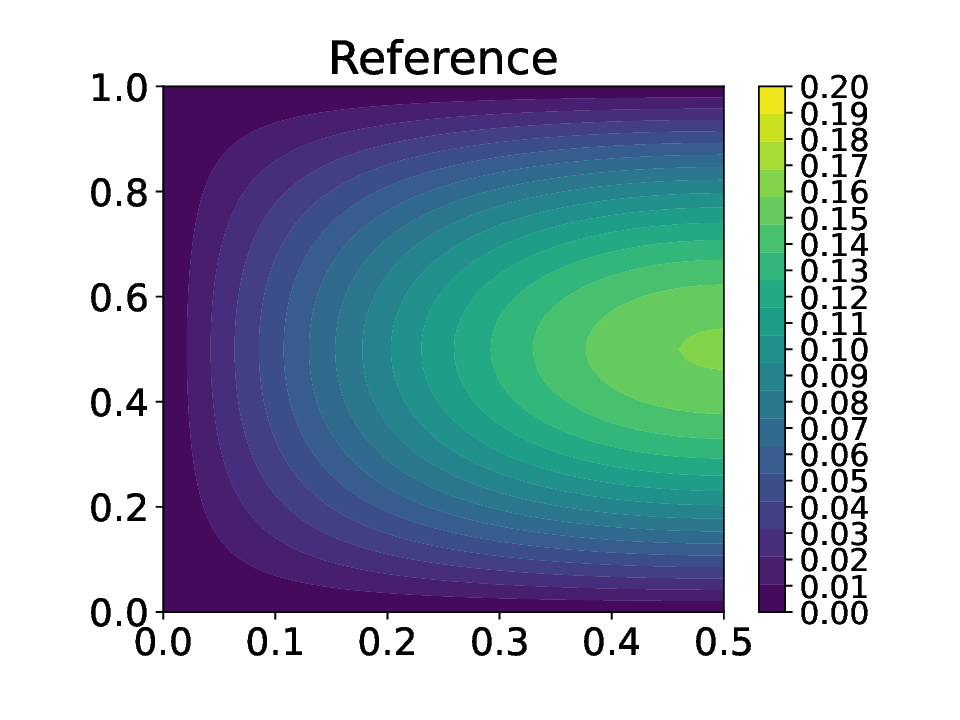 }\\
\scriptsize{(a) $D_1$}
\end{minipage}
\begin{minipage}{0.45\textwidth}
\centering
\includegraphics[width=6.5cm]{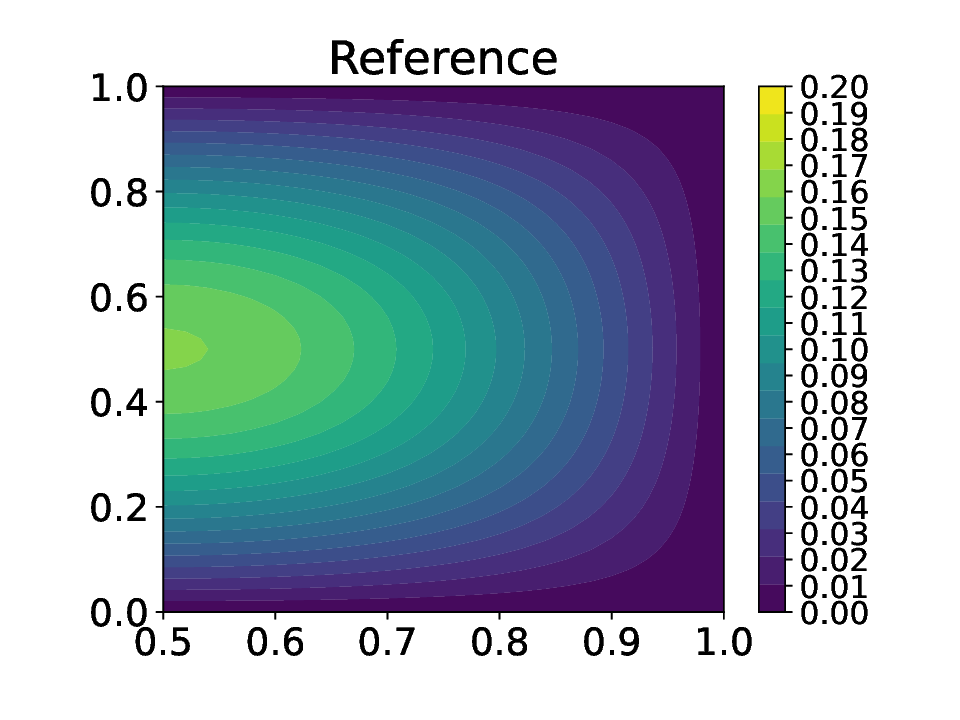 }\\
\scriptsize{(b) $D_2$}
\end{minipage}\\
\begin{minipage}{0.45\textwidth}
\centering
\includegraphics[width=6.5cm]{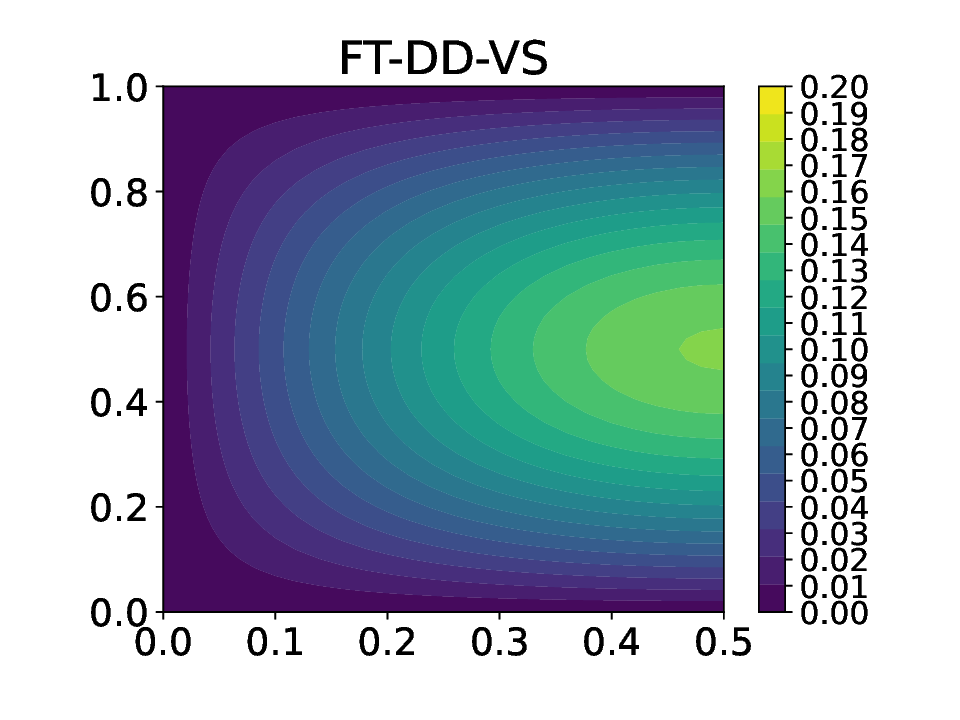 }\\
\scriptsize{(c) $D_1$}
\end{minipage}
\begin{minipage}{0.45\textwidth}
\centering
\includegraphics[width=6.5cm]{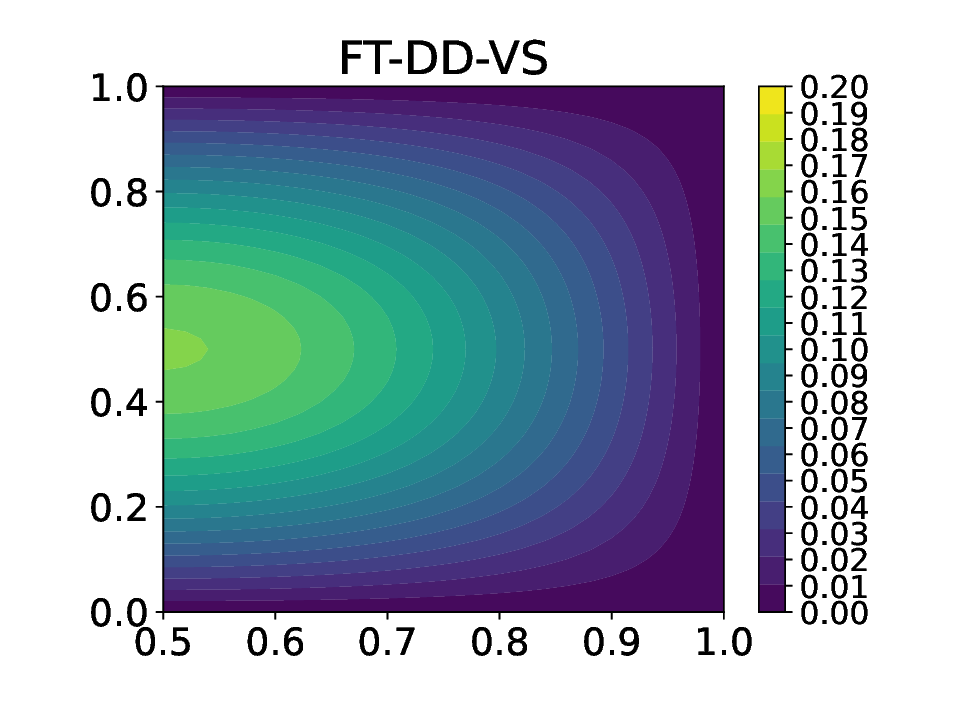 }\\
\scriptsize{(d) $D_2$}
\end{minipage}
\caption{Comparison of the solutions in subdomains for the FEM-BE and FT-DD-VS methods.}
\label{fig-heat-3}
\end{figure}

Fig.~\ref{fig-heat-3} presents the mean solutions at $t=1$, computed by the proposed method and the FEM-BE scheme in two subdomains. The top row shows the reference solutions, while the bottom row displays the corresponding results obtained using the proposed approach. It can be observed that the solution profiles and magnitudes produced by the proposed method closely match those of the reference.

Table \ref{table_heat} presents CPU time based on $M=10^3$ random samples, which includes the offline CPU time ($\mathcal{T}_{\text{off}}$), online CPU time ($\mathcal{T}_{\text{on}}$), total CPU time ($\mathcal{T}_{\text{tot}}$), and the average total CPU time ($T_{\text{tot}}$). The magnitude of the average total CPU time by the proposed method is much smaller than that of the FEM-BE method. It demonstrates that the proposed method can achieve a good trade-off in both approximation accuracy and computational efficiency for this numerical experiment.
\begin{table}[!tbh]
\caption{Comparison of the CPU time for FT-DD-VS and FEM-BE.}
\centering
\medskip\small\renewcommand{\arraystretch}{1.15}
\begin{tabular}{||ccccc||}
\hline
 Algorithm  & $\mathcal{T}_{\text{off}}$  &$\mathcal{T}_{\text{on}}$ & $\mathcal{T}_{\text{tot}}$ & $T_{\text{tot}}$ \\
\hline
FT-DD-VS  &   $2.90\times10^{2} s $ & $9.39 \times10^{2 }s$  &  $1.23 \times10^{3}s$ & $1.23s$ \\
FEM-BE & $\setminus $ & $\setminus$ & $3.98\times10^{5}s$  &$3.98\times10^{2}s$\\
\hline
\end{tabular}
\label{table_heat}
\end{table}

\subsection{Reaction diffusion equation}
In this subsection, we consider the reaction diffusion equation
\begin{equation}\label{eq:exa2}
   \frac{\partial u}{\partial t}(\bm{x},t;\bm{\xi})=\nabla\cdot (c_1(\bm{x};\bm{\xi}) \nabla u(\bm{x},t;\bm{\xi}))-c_2(\bm{x};\bm{\xi})u(\bm{x},t;\bm{\xi})+f(\bm{x},t;\bm{\xi}),
\end{equation}
with the initial and boundary conditions \eqref{eq:exa1}. The original domain $D=[0,1]\times [0,1]$ is divided into two subdomains $D_1 = [0,0.5]\times [0,1]$ and $D_2 = [0.5,1]\times [0,1]$.

\subsubsection{The first scenario}\label{exa:2}
First, we consider the source function 
\begin{equation*}
f(\bm{x},t;\bm{\xi})=e^{-t^2}e^{5(x_1+x_2)^2},    
\end{equation*}
and the coefficients
\begin{equation*}
c_1(\bm{x};\bm{\xi}) =\left\{
\begin{aligned}
100\xi_1, \ &\bm{x}\in D_1,\\
10\xi_2, \ &\bm{x}\in D_2,
\end{aligned}
\right.
\quad 
c_2(\bm{x};\bm{\xi}) = \left\{
\begin{aligned}
\xi_3, \ &\bm{x}\in D_1,\\
0.1\xi_4, \ &\bm{x}\in D_2,
\end{aligned}
\right.
\end{equation*}
where $\bm{\xi}:=({\xi}_1,{\xi}_2,{\xi}_3,{\xi}_4)\in[1,2]^{4}$.
The reference solution is calculated by FEM in space with the mesh size $h=0.05$ and the backward Euler scheme in time with the step size $\tau=5\times 10^{-4}$. We take $\omega^*=15$, $N_\omega = 15$, and  $|\Xi|=10$ for training in the offline stage of Algorithm \ref{algorithm-vs-sSAS}. In this problem, we select $N_{S_1} = N_{S_2} = N_{F_1} = N_{F_2} = 4$ and focus on the results of the interface problem and the subproblems.

\begin{figure}[!tbh]
\centering
\begin{minipage}{0.45\textwidth}
\centering
\includegraphics[width=6.5cm]{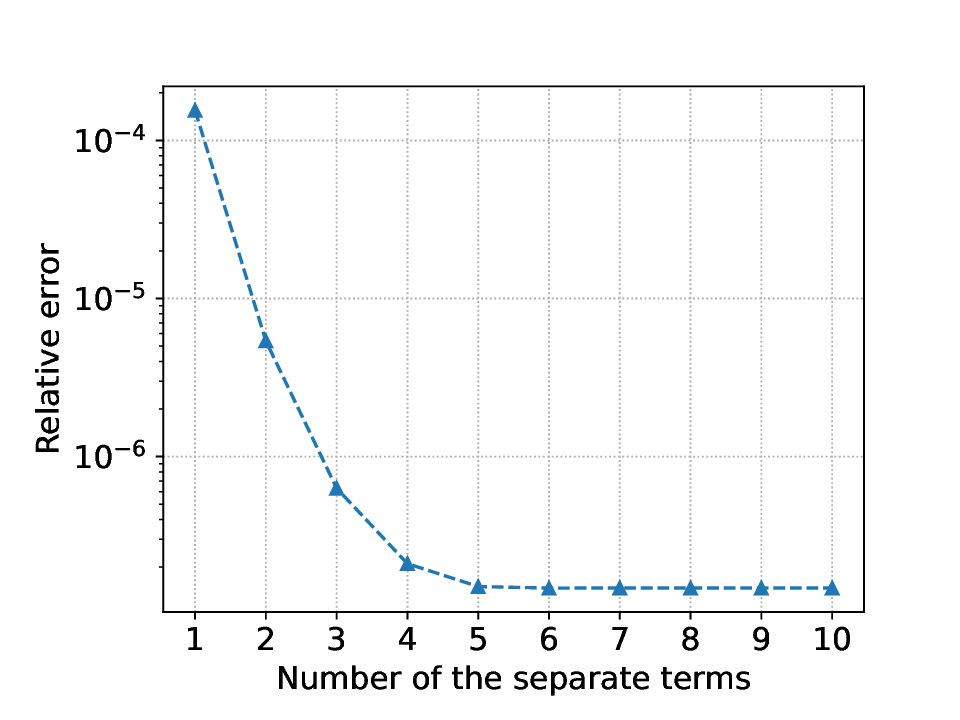 }\\
\scriptsize{(a) Interface problem}
\end{minipage}
\begin{minipage}{0.45\textwidth}
\centering
\includegraphics[width=6.5cm]{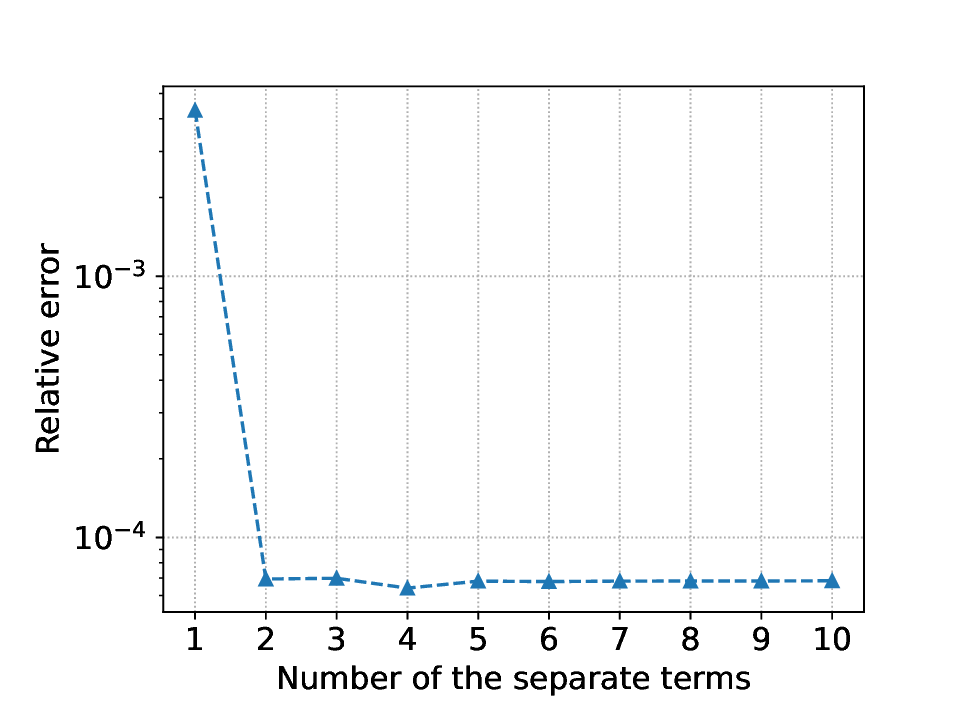 }\\
\scriptsize{(b) Subdomain problems}
\end{minipage}
\caption{Average relative error versus the number of separate terms.}
\label{fig2.1}
\end{figure}

Fig.~\ref{fig2.1} plots the average relative errors of the interface problem and the subproblems against the number of the separate terms $N$, where the average relative errors are computed based on $10^3$ samples. The relative error decreases as the number of separate terms increases, as shown in Fig.~\ref{fig2.1} (a) for the interface problem. For the subdomain problems, the relative error stabilizes when the number of separate terms is no less than 2.
From the figures, we can see that as the number of the separate terms increases, the approximation first becomes more accurate and then stable.

\begin{figure}[!tbh]
\centering
\begin{minipage}{0.45\textwidth}
\centering
\includegraphics[width=6.5cm]{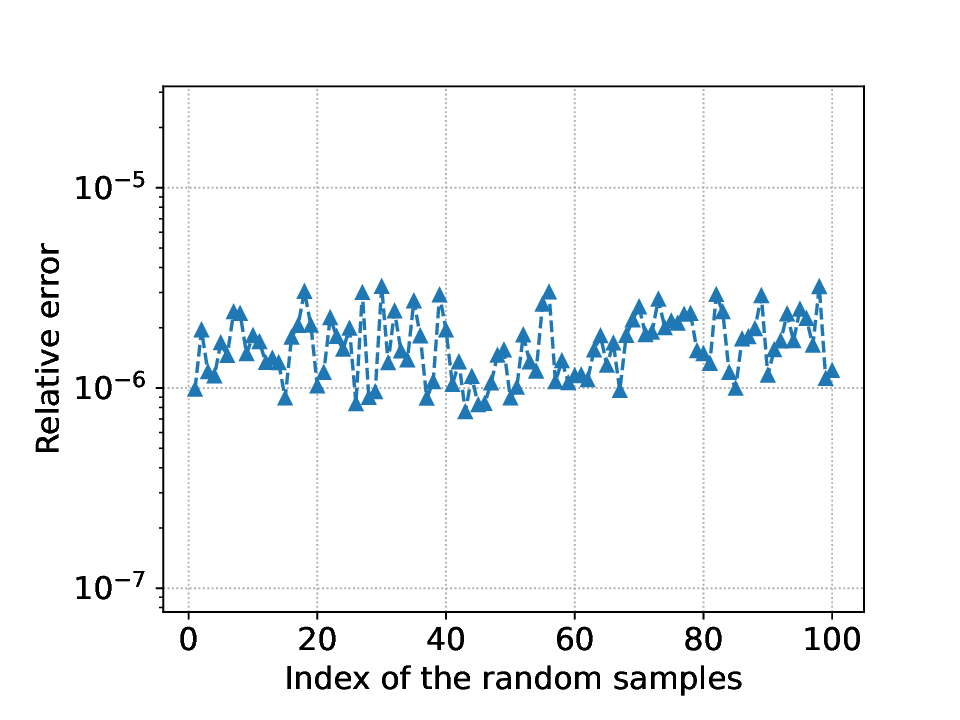 }\\
\scriptsize{(a) $t=0.2$}
\end{minipage}
\begin{minipage}{0.45\textwidth}
\centering
\includegraphics[width=6.5cm]{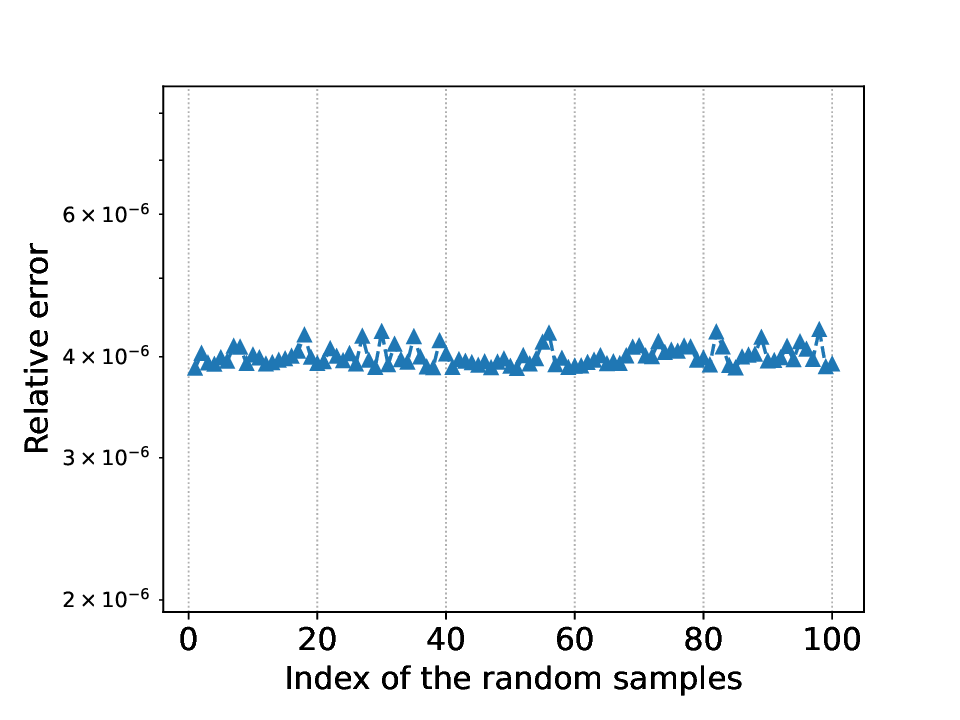 }\\
\scriptsize{(b) $t=0.8$}
\end{minipage}
\caption{Errors for the first $100$ samples.}
\label{fig2.2}
\end{figure}

\begin{figure}[!tb]
    \centering
    \includegraphics[width=7cm]{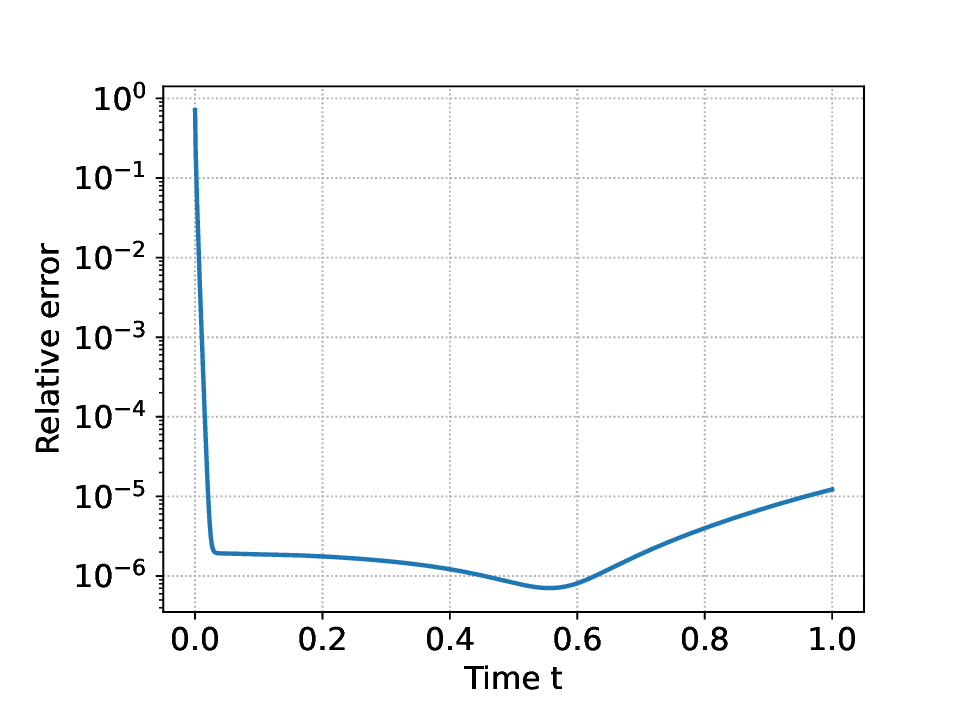 }
    \caption{Average relative errors versus time.}
    \label{fig2.3}
\end{figure}

To illustrate the individual relative error of the approximation in Fig. \ref{fig2.2}, we set $N_{\Gamma} = 3$ and $N_I=2$. Fig. \ref{fig2.2} demonstrates the relative errors of the first $100$ samples at $t=0.2$ and $t=0.8$. The range of the relative errors against the sample index at $t=0.2$ is wider than that at $t=0.8$.  Besides, we also plot the average relative errors versus time in Fig. \ref{fig2.3}. Here, the relative error first decreases and then increases one order of magnitude compared to that at $t=0.6$. This illustrates the proposed method can provide a good approximation for the original problem \eqref{eq:exa2} after $t=0.02$ with the relative error less than $10^{-5}$.

Finally, we present the CPU times based on $M=10^3$ random samples for both the proposed method and the FEM-BE method in Table \ref{table2}. The average total CPU time is reduced from $7.56\times 10^1s$ to $2.77\times 10^{-1}s$ using the proposed method.
In summary, the proposed method significantly reduces computational time compared to the FEM-BE method, while maintaining high accuracy.
\begin{table}[!tbh]
\caption{Comparison of the CPU time for FT-DD-VS and FEM-BE.}
\centering
\medskip\small\renewcommand{\arraystretch}{1.15}
\begin{tabular}{||ccccc||}
\hline
 Algorithm  & $\mathcal{T}_{\text{off}}$  &$\mathcal{T}_{\text{on}}$ & $\mathcal{T}_{\text{tot}}$ & $T_{\text{tot}}$ \\
\hline
FT-DD-VS  &   $2.87\times10^{1} s $ & $2.48 \times10^{2}s$  &  $2.77 \times10^{2}s$ & $2.77 \times10^{-1}s$  \\
\centering FEM-BE & $\setminus $ & $\setminus$ & $7.56\times10^{4}s$  &$7.56\times10^{1}s$\\
\hline
\end{tabular}
\label{table2}
\end{table}

\subsubsection{The second scenario}\label{exa:3}
In this example we set
\begin{equation*}
    f(\bm{x},t;\bm{\xi}) = \frac{1-t^2}{(1+t^2)^2}e^{5(x_1+x_2)^2}, \quad
c_1(\bm{x};\bm{\xi}) =\left\{
\begin{aligned}
     0,  &\ \bm{x}\in D_1, \\
    10\xi_1, &\ \bm{x}\in D_2,
\end{aligned}
 \right.
 \quad 
 c_2(\bm{x};\bm{\xi}) =\left\{
\begin{aligned}
  100\xi_2,  &\ \bm{x}\in D_1,\\
0,  &\ \bm{x}\in D_2,  
\end{aligned}
 \right.
\end{equation*}
and $\bm{\xi}:=({\xi}_1,{\xi}_2)\in[3,4]^{2} $.
The reference solution is calculated using FEM in space with a mesh size $h=0.05$ and the backward Euler scheme in time with a step size $\tau=10^{-3}$. In this problem, we take $\omega^*=15$, $N_\omega = 15$, $|\Xi|=10$, and set $N_{S_1} = N_{S_2} = N_{F_1} = N_{F_2} = 4$.

\begin{figure}[!tbh]
\centering
\begin{minipage}{0.45\textwidth}
\centering
\includegraphics[width=6.5cm]{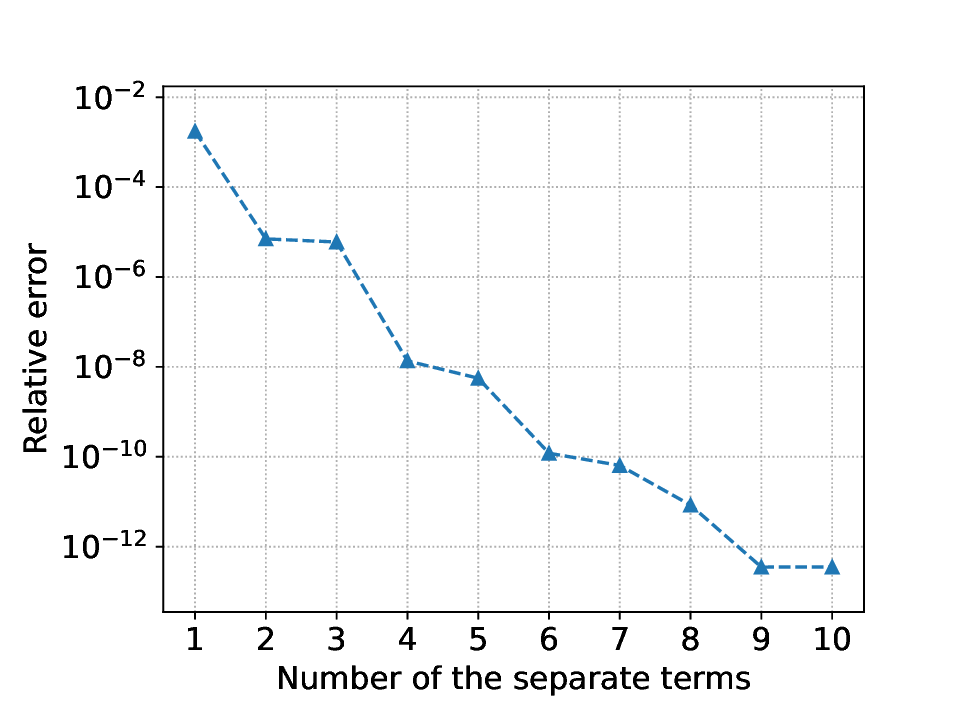 }\\
\scriptsize{(a) Interface problem}
\end{minipage}
\begin{minipage}{0.45\textwidth}
\centering
\includegraphics[width=6.5cm]{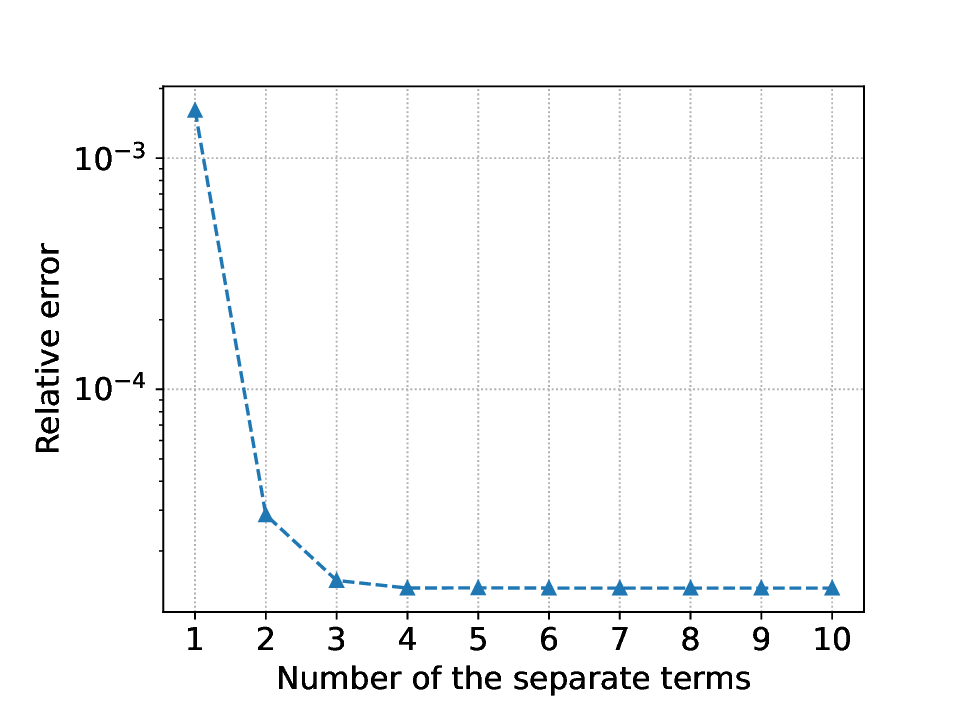 }\\
\scriptsize{(b) Subdomain problems}
\end{minipage}
\caption{Average relative error versus the number of separate terms.}
\label{fig5.1}
\end{figure}

Fig. \ref{fig5.1} presents the average relative errors of the interface problem and the subproblems against the number of the separate terms $N$. The average relative errors are calculated from $10^3$ random samples. The figures show that the approximation accuracy for the interface problem improves as $N$ increases. However, the error of the subproblems first decreases and then plateaus as $N$ increases.
Based on the results, we set $N_{\Gamma}=N_I=4$ to perform the Fourier inversion and ensure the approximation accuracy.

\begin{figure}[!tbh]
\centering
\begin{minipage}{0.45\textwidth}
\centering
\includegraphics[width=6.5cm]{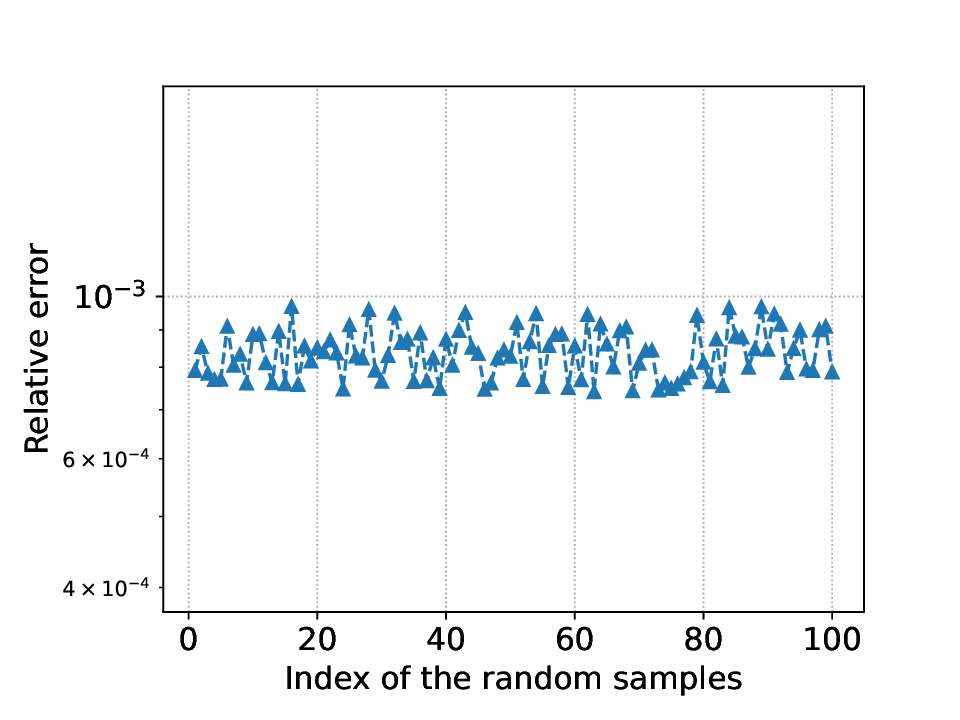 }\\
\scriptsize{(a) Errors for the first $100$ samples}
\end{minipage}
\begin{minipage}{0.45\textwidth}
\centering
\includegraphics[width=6.5cm]{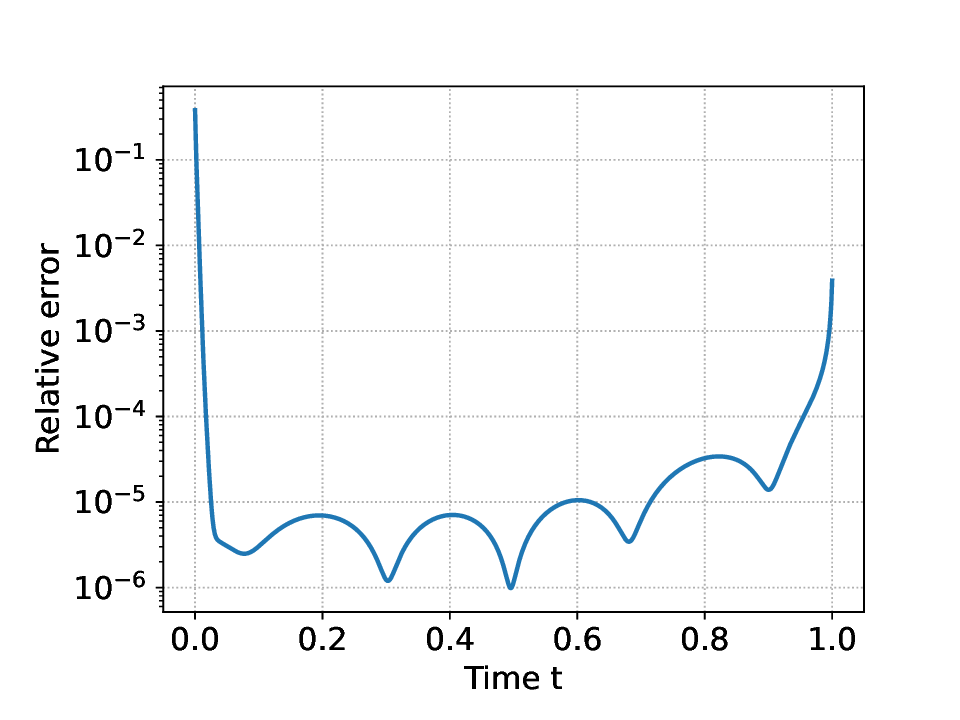 }\\
\scriptsize{(b) Average relative errors versus the time}
\end{minipage}
 \caption{Comparison of the errors.}
    \label{fig5.2}
\end{figure}

To illustrate the approximation accuracy, the relative errors of the first 100 samples along with the errors versus the time are plotted in Fig.~\ref{fig5.2}. The results in Fig.~\ref{fig5.2} (a) demonstrate that the proposed method provides accurate approximations for the first 100 random samples, because the curve of relative errors stay below the threshold of $10^{-3}$.
In Fig.~\ref{fig5.2} (b), we can see that the relative errors within $t\in[0.1, 0.9]$ exhibits a clear periodic variation pattern. Moreover, the average relative errors lies in the interval $[10^{-6}, 10^{-4}]$ when $t\in[0.02, 0.96]$. In this example, the proposed method shows a good performance for $t\in[0.02, 0.96]$.

Finally, Table \ref{table3} compares the proposed method with the FEM-BE approach in terms of the computational cost. We find that the CPU time of the proposed method is orders of magnitude smaller than that of the FEM-BE method. The results demonstrate that the proposed method maintains high accuracy while achieving significant computational time reduction compared to the FEM-BE approach.
\begin{table}[!tbh]
\caption{Comparison of the CPU time for FT-DD-VS and FEM-BE.}
\centering
\medskip\small\renewcommand{\arraystretch}{1.15}
\begin{tabular}{||ccccc||}
\hline
 Algorithm  & $\mathcal{T}_{\text{off}}$  &$\mathcal{T}_{\text{on}}$ & $\mathcal{T}_{\text{tot}}$ & $T_{\text{tot}}$ \\
\hline
 FT-DD-VS  &   $2.81\times10^{1} s $ & $2.02 \times10^{2}s$  &  $2.30 \times10^{2}s$ & $2.30 \times10^{-1}s$  \\
 FEM-BE & $\setminus $ & $\setminus$ & $3.69\times10^{4}s$  &$3.69\times10^{1}s$\\
\hline
\end{tabular}
\label{table3}
\end{table}

\section{Conclusions}
\label{sec-Conclusions}
We propose a model order reduction method for parametric dynamical systems, in which time-dependent problems are first transformed into frequency-variable elliptic equations using the Fourier transform. The converted equations become time-independent, allowing for parallel computation across different frequency values. To construct an efficient and accurate separable approximation for the resulting complex-valued elliptic equations, we employ the variable-separation-based domain decomposition method.
This approach leverages domain decomposition to convert the solution of the complex elliptic problem into solving interface problem and subproblems, making it well-suited for handling complex scenarios. Subsequently, the VS method is applied to establish surrogate models for both the interface problem and the subproblems. Based on this framework, the solution to the time-dependent problem is obtained by applying the inverse Fourier transform to the elliptic equation solutions.
A key advantage of the method is its ability to decouple offline and online computations, ensuring that the online phase remains entirely independent of spatial discretization, thus achieving high efficiency. As a result, the overall computational cost of the proposed method is significantly lower than that of standard finite element method, a fact supported by numerical experiments.
Nevertheless, the Fourier transform has inherent limitations, since the functions to be transformed should satisfy certain assumptions.
Consequently, this method faces significant constraints when applied to time-dependent problems. In future work, we aim to explore domain decomposition directly on time-dependent problems to achieve efficient computation, and focus on extending to general parameter distributions by adapting the subdomain partitioning strategy.

\section*{Acknowledgments}
The authors would like to thank the two anonymous referees for their helpful comments that improved the quality of the manuscript.

This work was supported by the National Key R \& D Program of China (2021YFA001300),
the National Natural Science Foundation of China (12401567, 12271150, 12471405),
the Hunan Provincial Natural Science Foundation of China (2023JJ10001),
the Science and Technology Innovation Program of Hunan Province (2022RC1190), and
the Hunan Provincial Innovation Foundation for Postgraduate (CX20250499).

\end{document}